\documentclass[reqno]{amsart}
\usepackage{mathrsfs}
\usepackage{amssymb}
\usepackage{latexsym}
\usepackage{yfonts}
\usepackage{mathrsfs}
\usepackage{natbib}
\usepackage{amsmath}
\usepackage{graphicx,color}
\usepackage[pdftex,plainpages=false,colorlinks,hyperindex,bookmarksopen,linkcolor=blue,citecolor=red,urlcolor=red]{hyperref}
\usepackage[final]{showkeys}
\usepackage{bm}
\usepackage{marvosym}
\DeclareMathAlphabet{\mathpzc}{OT1}{pzc}{m}{it}
\usepackage{bm}

\bibpunct{[}{]}{;}{n}{,}{,}

\newtheorem{te}{Theorem}[section]

\newtheorem{os}{Remark}[section]

\numberwithin{equation}{section}
\numberwithin{te}{section}

\allowdisplaybreaks

\def \l {\left(}
\def \r {\right)}
\def \ll {\left\lbrace}
\def \rr {\right\rbrace}
\def \lll {\left|}
\def \rrr {\right|}

\begin{document}

	\title [Space-time fractional telegraph equations] {Time-changed processes governed by space-time fractional telegraph equations}
	\author{Mirko D'Ovidio}
	\address{Department of basic and applied sciences for engineers, Sapienza University of Rome}
	\email {mirko.dovidio@uniroma1}
	\author{Enzo Orsingher} 
	\address{Department of Statistical Sciences, Sapienza University of Rome}
		\email {enzo.orsingher@uniroma1.it}
	\author{Bruno Toaldo}
	\address{Department of Statistical Sciences, Sapienza University of Rome}
	\email {bruno.toaldo@uniroma1.it}
	 
	\date{\today}

		\begin{abstract}
In this work we construct compositions of vector processes of the form $\bm{S}_n^{2\beta} \l c^2 \mathpzc{L}^\nu (t) \r$, $t>0$, $\nu \in \l 0, \frac{1}{2} \right]$, $\beta \in \l 0,1 \right]$, $n \in \mathbb{N}$, whose distribution is related to space-time fractional $n$-dimensional telegraph equations. We present within a unifying framework the pde connections of $n$-dimensional isotropic stable processes $\bm{S}_n^{2\beta}$ whose random time is represented by the inverse $\mathpzc{L}^\nu (t)$, $t>0$, of the superposition of independent positively-skewed stable processes, $\mathpzc{H}^\nu (t) = H_1^{2\nu} (t) + \l 2\lambda \r^{\frac{1}{\nu}} H_2^\nu (t)$, $t>0$, ($H_1^{2\nu}$, $H_2^\nu$, independent stable subordinators). As special cases for $n=1$, $\nu = \frac{1}{2}$ and $\beta = 1$ we examine the telegraph process $T$ at Brownian time $|B|$ (\citet{ptrf}) and establish the equality in distribution $B \l c^2 \mathpzc{L}^{\frac{1}{2}} (t) \r \, \stackrel{\textrm{law}}{=} \, T \l \lll B(t) \rrr \r$, $t>0$. Furthermore the iterated Brownian motion (\citet{allouba}) and the two-dimensional motion at finite velocity with a random time are investigated. For all these processes we present their counterparts as Brownian motion at delayed stable-distributed time.
		\end{abstract}

\maketitle

\tableofcontents
	
	\section{Introduction and preliminaries}

	\subsection{Introduction}

The study of the interplay between fractional equations and stochastic processes has began in the middle of the Eighties with the analysis of simple time-fractional diffusion equations (see \citet{fujita} for a rigorous work on this field, or more recently \citet{nane}, where the compositions of Brownian sheets with Brownian motions are considered). In some papers the connection between fractional diffusion equations and stable processes is explored (see, for example, \citet{orsann, Zol86}). The iterated Brownian motion has distribution satisfying the following fractional equation
\begin{equation} 
	\frac{\partial^{\frac{1}{2}}}{\partial t^{\frac{1}{2}}} u(x, t) \, = \, \frac{1}{2^{\frac{3}{2}}} \frac{\partial^2}{\partial x^2} u(x, t), \qquad x \in \mathbb{R}, t>0,
	\label{fractionaldiffusion}
\end{equation}
(see for example \citet{allouba}) and also the fourth-order equation
\begin{equation} 
	\frac{\partial}{\partial t} u(x, t) \, = \, \frac{1}{2^3} \frac{\partial^4}{\partial x^4} u(x, t)+ \frac{1}{2\sqrt{2\pi t}} \frac{d^2}{dx^2} \delta (x), \qquad x \in \mathbb{R}, t>0,
	\label{quartoordine}
\end{equation}
see \citet{debl} (also for an interpretation of the iterated Brownian motion to model the motion of a gas in a crack).

When the fractional equation has a telegraph structure, with more than one time-fractional derivative involved, that is for $\nu \in \l 0, 1\right]$
\begin{equation} 
	\l \frac{\partial^{2\nu}}{\partial t^{2\nu}} + 2\lambda \frac{\partial^\nu}{\partial t^\nu} \r u(x, t) \, = \, c^2 \frac{\partial^2}{\partial x^2} u(x, t), \qquad x \in \mathbb{R}, t>0, \lambda >0, c >0,
	\label{}
\end{equation}
the relationship of its solution with the time-changed telegraph processes is examined and established in \citet{ptrf} . The space-fractional telegraph equation (with M. Riesz space derivatives) has been considered in \citet{zhao}, while the connection between space-fractional equations and asymmetric stable processes has been established in \citet{Feller52}.

Fractional telegraph equations from the analytic point of view have been studied by many authors (see \citet{saxena} for equations with $n$ time derivatives). For their solutions have been worked out also numerical techniques (see, for example, \citet{moma}).
Telegraph equations have an extraordinary importance in electrodynamics (the scalar Maxwell equations are of this type), in the theory of damped vibrations and in probability because they are connected with finite velocity random motions.

In this paper we consider various types of processes obtained by composing symmetric stable processes $\bm{S}_n^{2\beta} (t)$, $t>0$, $0<\beta \leq 1$, with the inverse of the sum of two independent stable subordinators (instead of one as in \citet{baem}) say $\mathpzc{L}^\nu (t)$, $t>0$, $0< \nu \leq \frac{1}{2}$. These time-changed processes, $\bm{W}_n (t) = \bm{S}_n^{2\beta} \l c^2 \mathpzc{L}^\nu (t) \r$, $t>0$, have distributions, $w_\nu^\beta (\bm{x}, t)$, $\bm{x} \in \mathbb{R}^n$, $t>0$, which satisfy telegraph-type space-time fractional equations of the form
\begin{equation} 
	\l \frac{^C\partial^{2\nu}}{\partial t^{2\nu}} + 2\lambda \frac{^C\partial^\nu}{\partial t^\nu} \r w_\nu^\beta \l \bm{x}, t \r \, = \, -c^2 \l -\Delta \r^\beta w_\nu^\beta \l \bm{x}, t\r, \qquad \bm{x} \in \mathbb{R}^n, t>0, c>0, \lambda >0,
	\label{1.1}
\end{equation}
where $0 < \beta \leq 1$, $0<\nu \leq \frac{1}{2}$, subject to the initial condition
\begin{equation} 
	w_\nu^\beta \l \bm{x}, 0 \r \, = \, \delta (\bm{x}).
	\label{initial condition}
\end{equation}
The fractional Laplacian $\l -\Delta \r^\beta$, appearing in \eqref{1.1}, is defined and analyzed in Section \ref{laplaciansec} below. The fractional derivatives appearing in \eqref{1.1} are meant in the Dzerbayshan-Caputo sense, that is, for an absolutely continuous function $f \in L^1 \l \mathbb{R} \r$ (for fractional calculus consult \citet{kill}),
\begin{equation} 
	\frac{^C\partial^\nu}{\partial t^\nu} f(t) \, = \, \frac{1}{\Gamma \l m-\nu \r} \int_0^t \frac{\frac{d^m}{ds^m}f(s)}{(t-s)^{\nu + 1-m}} \, ds, \qquad m-1 < \nu < m, m \in \mathbb{N}.
	\label{}
\end{equation}
Equation \eqref{1.1} includes as particular cases all fractional equations studied so far (including diffusion equations) and also the main equations of mathematical physics as limit cases. Thus the distribution of the composed process $\bm{S}_n^{2\beta} \l \mathpzc{L}^\nu (t) \r$, $t>0$, represents the fundamental solution of the most general $n$-dimensional time-space fractional telegraph equation. We give the general Fourier transform of the solution to \eqref{1.1} with initial condition \eqref{initial condition} as
\begin{align}
	&\mathbb{E}e^{i \bm{\xi} \cdot \bm{S}_n^{2\beta} \l c^2 \mathpzc{L}^\nu (t) \r} \, = \notag \\
	 = \, & \frac{1}{2} \left[ \l 1+ \frac{\lambda}{\sqrt{\lambda^2- c^2 \left\| \bm{\xi} \right\|^{2\beta}}} \r E_{\nu,1} \l r_1 t^\nu \r + \l 1-\frac{\lambda}{\sqrt{\lambda^2 - c^2 \left\| \bm{\xi} \right\|^{2\beta}}} \r E_{\nu,1} \l r_2 t^\nu \r  \right],
	\label{}
\end{align}
where
\begin{equation} 
	r_1 \, = \, -\lambda + \sqrt{\lambda^2 - c^2 \left\| \bm{\xi} \right\|^{2\beta}}, \qquad r_2 \, = \, -\lambda -\sqrt{\lambda^2 - c^2 \left\| \bm{\xi} \right\|^{2\beta}}.
	\label{}
\end{equation}
and
\begin{equation} 
	E_{\nu, \psi} \l x \r \, = \, \sum_{k=0}^\infty \frac{x^k}{\Gamma \l \nu k + \psi \r}, \qquad \nu, \psi > 0,
	\label{def mittag leffler}
\end{equation}
is the two-parameters Mittag-Leffler function (see, for example, \citet{hms} for a general overview on the Mittag-Leffler functions). Our result therefore includes all previous results in a unique framework and sheds an additional insight into the literature in this field.

An important role in our analysis is played by the time change based on the process $\mathpzc{L}^\nu (t)$, $t>0$.
We consider first the sum of two independent positively skewed stable r.v.'s $H_1^{2\nu} (t)$ and $H_2^\nu (t)$, $t>0$, $0< \nu \leq \frac{1}{2}$,
\begin{equation} 
	\mathpzc{H}^\nu (t) \, = \, H_1^{2\nu} (t) + \l 2\lambda \r^{\frac{1}{\nu}} H_2^\nu (t), \qquad t>0,
	\label{}
\end{equation}
whose distribution $\mathpzc{h}_\nu (x, t)$ is governed by the space fractional equation
\begin{equation} 
	\frac{\partial}{\partial t} \mathpzc{h}_\nu (x, t) \, = \, - \l \frac{\partial^{2\nu}}{\partial x^{2\nu}} + 2\lambda \frac{\partial^\nu}{\partial x^\nu} \r \, \mathpzc{h}_\nu (x, t), \qquad x \geq 0, t>0, 0<\nu \leq \frac{1}{2}.
	\label{pallino cerchiato}
\end{equation}
In \eqref{pallino cerchiato} the fractional derivatives must be meant in the Riemann-Liouville sense which, for a function $f \in L^1 \l \mathbb{R} \r$, is defined as
\begin{equation} 
	\frac{\partial^\nu}{\partial x^\nu} f(x) \, = \, \frac{1}{\Gamma \l m-\nu \r} \frac{d^m}{dx^m} \int_0^x \frac{f(s)}{(x-s)^{\nu + 1 - m}} \, ds, \qquad m-1 < \nu < m, m \in \mathbb{N}.
	\label{}
\end{equation}
We then take the inverse $\mathpzc{L}^\nu (t)$, $t>0$, to the process $\mathpzc{H}^\nu (t)$, $t>0$, defined as
\begin{equation} 
	\mathpzc{L}^\nu (t) \, = \, \inf \ll s>0: H_1^{2\nu} (s) + \l 2\lambda \r^{\frac{1}{\nu}} H_2^\nu (s) \geq t \rr, \qquad t>0,
	\label{119}
\end{equation}
whose distribution is related to that of  $\mathpzc{H}^\nu (t)$, $t>0$, by means of the formula
\begin{equation} 
	\Pr \ll \mathpzc{L}^\nu (t) < x \rr \, = \, \Pr \ll \mathpzc{H}^\nu (x) > t \rr.
	\label{}
\end{equation}
The distribution $\mathpzc{l}_\nu (x, t)$ of $\mathpzc{L}^\nu (t)$, $t>0$, satisfies the time-fractional telegraph equation
\begin{equation} 
	\l \frac{\partial^{2\nu}}{\partial t^{2\nu}} + 2\lambda \frac{\partial^\nu}{\partial t^\nu} \r \mathpzc{l}_\nu (x, t) \, = \, -\frac{\partial}{\partial x} \mathpzc{l}_\nu (x, t), \qquad x \geq 0, t>0, 0 < \nu \leq \frac{1}{2},
	\label{triangolo}
\end{equation}
where the fractional derivatives appearing in \eqref{triangolo} are again in the Riemann-Liouville sense.
We are able to give explicit forms of the Laplace transforms of $\mathpzc{h}_\nu (x, t)$ and $\mathpzc{l}_\nu (x, t)$ in terms of Mittag-Leffler functions for all values of $0 < \nu \leq \frac{1}{2}$. For example, for the distribution $\mathpzc{l}_\nu (x, t)$ of $\mathpzc{L}^\nu (t)$ we have that, for $\gamma < \lambda^2$,
\begin{align}
	&  \int_0^\infty e^{-\gamma x} \mathpzc{l}_\nu (x, t) \, dx \, = \notag \\
	 = \, & \frac{1}{2} \left[ \l 1+ \frac{\lambda}{\sqrt{\lambda^2- \gamma}} \r E_{\nu,1} \l r_1 t^\nu \r + \l 1-\frac{\lambda}{\sqrt{\lambda^2 - \gamma}} \r E_{\nu,1} \l r_2 t^\nu \r  \right],
	\label{}
\end{align}
where 
\begin{equation} 
	r_1 \, = \, -\lambda + \sqrt{\lambda^2 - \gamma}, \qquad r_2 \, = \, -\lambda -\sqrt{\lambda^2 - \gamma}.
	\label{}
\end{equation}
The distribution $\mathpzc{l}_\nu (x, t)$ of $\mathpzc{L}^\nu (t)$, $t>0$, has the general form
\begin{equation} 
	\mathpzc{l}_\nu (x, t) \, = \, \int_0^t l_{2\nu} \l  x, s \r \, h_\nu (t-s, 2\lambda x) \, ds + 2\lambda \int_0^t l_\nu (2\lambda x, s) \, h_{2\nu} \l t-s, x \r \, ds,
	\label{}
\end{equation}
where the distributions of $H^{2\nu}$, $H^\nu$, and that of their inverse processes $L^{2\nu}$ and $L^\nu$ appear. For our analysis it is relevant to obtain the distributions of $\mathpzc{H}^{\frac{1}{2}} (t)$, $t>0$, and $\mathpzc{L}^{\frac{1}{2}} (t)$, $t>0$. We also obtain explicitely the distributions of $H^{\frac{1}{3}} (t)$ and $H^{\frac{2}{3}} (t)$, $t>0$, and also of their inverses $L^{\frac{1}{3}} (t)$ and $L^{\frac{2}{3}} (t)$, $t>0$, in terms of Airy functions. By means of the convolutions of these distributions we arrive at the following cumbersome density of the random time $\mathpzc{L}^{\frac{1}{3}} (t)$, $t>0$,
\begin{align} 
	\Pr \ll \mathpzc{L}^{\frac{1}{3}} (t) \in dx \rr \,  = \, &  \frac{2\lambda}{\sqrt{\pi}} \int_0^t ds \int_0^\infty dw \, e^{-w} w^{-\frac{1}{6}} \, \textrm{Ai} \l -x \sqrt[3]{\frac{2^2w}{3(t-s)^2}} \r \, \textrm{Ai} \l \frac{2\lambda x}{\sqrt[3]{3s}} \r \bm{\cdot} \notag \\
	&  \bm{\cdot} \frac{3}{\sqrt[3]{3s}} \sqrt[3]{\frac{2^2}{3(t-s)^2}} \left[ \frac{x}{2s} + \frac{s}{t-s} \right] \; dx.
	\label{l di un terzo intro}
\end{align}

For $n=1$, $\beta =1$ and $\nu =1$ in \eqref{1.1}, we get the telegraph equation which is satisfied by the distribution of the one-dimensional telegraph process
\begin{equation} 
	T(t) \, = \, V(0) \int_0^t (-1)^{N(s)} \, ds, \qquad t>0,
	\label{costruzionetelegrafo}
\end{equation}
where $N(t)$, $t>0$ is an homogeneous Poisson process, with parameter $\lambda >0$, independent from the symmetric r.v. $V(0)$ (with values $\pm c$). Properties of this process (including first-passage time distributions) are studied in \citet{kanno} and a telegraph process with random velocities has been recently considered by \citet{stad}.

For $n=1$, $\beta = 1$ and $\nu = \frac{1}{2}$ the special equation
\begin{align}
	\begin{cases}
	\l \frac{\partial}{\partial t} + 2\lambda \frac{\partial^{\frac{1}{2}}}{\partial t^{\frac{1}{2}}} \r w_{\frac{1}{2}}^1 (x, t) \, = \, c^2 \frac{\partial^2}{\partial x^2} w_{\frac{1}{2}}^1 (x, t), \qquad x \in \mathbb{R}, t>0, \\
	w_{\frac{1}{2}}^1 (x, 0) \, = \, \delta (x),
	\end{cases}
	\label{orsbegh}
\end{align}
has solution coinciding with the distribution of $T \l \left| B(t) \right| \r$, $t>0$, where $ \left| B(t) \right|$, $t>0$, is a reflecting Brownian motion independent from $T$ (see \citet{ptrf}). For $\lambda \to \infty$, $c \to \infty$, in such a way that $\frac{c^2}{\lambda} \to 1$ the fractional diffusion equation \eqref{fractionaldiffusion} is obtained from \eqref{orsbegh} and the composition $T \l \left| B(t) \right| \r$, $t>0$, converges in distribution to the iterated Brownian motion. Our result, specialized to this particular case gives the following unexpected equality in distribution
\begin{equation} 
	T \l \left| B(t) \right| \r \, \stackrel{\textrm{law}}{=} \, B \l c^2 \mathpzc{L}^{\frac{1}{2}} (t) \r, \qquad t>0,
	\label{roletelegraph}
\end{equation}
where
\begin{equation} 
	\Pr \ll B \l c^2 \mathpzc{L}^\nu (t) \r \in dx \rr \, = \, \frac{\lambda \, dx}{c \pi} \int_0^t \frac{1}{\sqrt{s(t-s)}} e^{-\frac{x^2}{4c^2s} - \frac{\lambda^2 s^2}{t-s}} \l \frac{s}{2(t-s)} +1 \r \, ds,
	\label{}
\end{equation}
and
\begin{equation} 
	\Pr \ll T \l \left| B(t) \right| \r \in dx \rr \, = \, \int_0^\infty \Pr \ll T(s) \in dx \rr \, \Pr \ll \lll B(t) \rrr \in ds \rr.
	\label{}
\end{equation}
The absolutely continuous component of the distribution of the telegraph process $T(t)$, $t>0$, reads
\begin{equation} 
	\Pr \ll T(s) \in dx \rr \, = \, \frac{dx \, e^{-\lambda t}}{2c} \ll \lambda \, I_0 \l \frac{\lambda}{c} \sqrt{c^2t^2 - x^2} \r + \frac{\partial}{\partial t} \, I_0 \l \frac{\lambda}{c} \sqrt{c^2t^2-x^2} \r  \rr,  
	\label{}
\end{equation}
where $|x| < ct$, $t>0$, $c>0$, and
\begin{equation} 
	I_0 (x) \, = \, \sum_{k=0}^\infty \l \frac{x}{2} \r^{2k} \frac{1}{(k!)^2}.
	\label{}
\end{equation}

For $n=2$, $\beta = 1$ and $\nu = 1$, equation \eqref{1.1} coincides with that of damped planar vibrations (we call it planar telegraph equation) and governs the vertical oscillations of thin deformable structures. The solution to
\begin{align}
	\begin{cases}
	\l \frac{\partial^2}{\partial t^2} + 2\lambda \frac{\partial}{\partial t} \r r (x, y, t) \, = \, c^2 \l \frac{\partial^2}{\partial x^2} + \frac{\partial^2}{\partial y^2} \r r(x, y, t), \qquad  x^2+y^2 < c^2t^2, t>0, \\
	r(x, y, 0) \, = \, \delta(x, y), \\
	r_t (x, y, 0) \, = \, 0,
	\end{cases}
	\label{2}
\end{align}
corresponds to the distribution $r(x, y, t)$ of the vector $\bm{T} (t) = \l X(t), Y(t) \r$ related to a planar motion described in \citet{orsdegr}.
This random motion $\bm{T}(t)$, $t>0$, is performed at finite velocity $c$, possesses sample paths composed by segments whose orientation is uniform in $\l 0, 2\pi \r$, and with changes of direction at Poisson times. The distribution $r(x, y, t)$ of $\bm{T}(t)$, $t>0$, is concentrated inside a circle $C_{ct}$ of radius $ct$ and has an absolutely continuous component which reads
\begin{equation} 
	r(x, y , t) \, = \, \frac{\lambda}{2 \pi c} \frac{e^{-\lambda t + \frac{\lambda}{c} \sqrt{c^2t^2 - \l x^2+y^2 \r}}}{\sqrt{c^2t^2 - \l x^2+y^2 \r}}, \qquad (x,y) \in C_{ct}, t>0.
	\label{}
\end{equation}
If no Poisson event occurs, the moving particle reaches the boundary $\partial C_{ct}$ of $C_{ct}$ with probability $e^{-\lambda t}$.
The vector process $\bm{T}(t)$, $t>0$, taken at a random time represented by a reflecting Brownian motion, $\left| B(t) \right|$, has distribution
\begin{equation} 
	q(x, y, t) \, = \, \int_0^\infty \Pr \ll X(t) \in ds, Y(t) \in ds \rr \, \Pr \ll \left| B(t) \right| \in ds \rr
	\label{4}
\end{equation}
which satisfies the fractional equation
\begin{equation} 
	\l \frac{\partial}{\partial t} + 2\lambda \frac{\partial^{\frac{1}{2}}}{\partial t^{\frac{1}{2}}} \r q(x, y, t) \, = \, c^2 \l \frac{\partial^2}{\partial x^2} + \frac{\partial^2}{\partial y^2} \r q(x, y, t), \qquad (x,y) \in \mathbb{R}^2, t>0.
	\label{}
\end{equation}
However, the distribution of $\bm{B}_2 \l c^2 \mathpzc{L}^{\frac{1}{2}} (t) \r$, $t>0$, does not coincide with \eqref{4} ($\bm{B}_2$ is a two dimensional Brownian motion). In this case the role of $T(t)$, $t>0$, in \eqref{roletelegraph} is here played by a process which is a slight modification of $\bm{T} (t)$, $t>0$. We take the planar process with law
\begin{equation} 
	\mathfrak{r}(x, y, t) \, = \, \frac{\lambda \, e^{-\lambda t}}{2\pi c} \left[ \frac{e^{\frac{\lambda}{c} \sqrt{c^2t^2 - \l x^2+y^2 \r}} + e^{-\frac{\lambda}{c} \sqrt{c^2t^2 - \l x^2+y^2 \r}}}{\sqrt{c^2t^2 - \l x^2+y^2 \r}} \right], \quad x^2+y^2 < c^2t^2, t>0,
	\label{6}
\end{equation}
which also solves equation \eqref{2}. The process with distribution
\begin{align}
	\mathfrak{q}(x, y, t) \, & = \, \int_0^\infty \mathfrak{r}(x, y, s) \left[ \Pr \ll \lll B(t) \rrr \in ds \rr + \frac{1}{2\lambda} \frac{\partial^{\frac{1}{2}}}{\partial t^{\frac{1}{2}}} \Pr \ll \lll B(t) \rrr \in ds \rr \right] \notag \\
	& = \, \int_0^\infty \l \mathfrak{r} (x, y, s) + \frac{\partial}{\partial s} \mathfrak{r} (x, y, s) \r \Pr \ll \lll B(t) \rrr \in ds  \rr,
	\label{}
\end{align}
has the same law of a planar Brownian motion at the time $\mathpzc{L}^{\frac{1}{2}} (t)$, $t>0$.
The process $\mathfrak{T} (t)$, $t>0$, possessing distribution \eqref{6} is obtained from $\bm{T}(t)$, $t>0$, by disregarding displacements started off by even-order Poisson events.

\subsection{Notations}
For the reader convenience we list below the main notations used throughout the paper.
	\begin{enumerate}
	\item[$\bullet$] $\bm{S}_n^{2\beta} (t ) = \l S_1^{2\beta} (t ), S_2^{2\beta} (t), \cdots, S_n^{2\beta} (t) \r$, $t>0$, $0<\beta \leq 1$, $ n \in \mathbb{N}$  is a isotropic stable $n$-dimensional process with law $v_\beta \l \bm{x}, t \r$, $\bm{x} \in \mathbb{R}^n$, $t>0$.
	\item[$\bullet$] $H^\nu (t)$, $t>0$, $0< \nu < 1$,  is a totally positively-skewed stable process (stable subordinator), with law $h_\nu (x, t)$, $x\geq 0$, $t>0$.
		\item[$\bullet$] $L^\nu (t)$, $t>0$, is the inverse of $H^\nu (t)$, $t>0$, and has law $l_\nu (x, t)$, $x\geq 0$, $t>0$.
	\item[$\bullet$] $\mathpzc{H}^\nu (t)  = H_1^{2\nu} (t) + \l 2\lambda \r^{ \frac{1}{\nu} } H_2^\nu (t) $, $t>0$,  is the sum of two independent stable subordinators and has law $\mathpzc{h}_\nu (x, t)$, $x \geq 0$, $t>0$.
	\item[$\bullet$] $\mathpzc{L}^\nu (t)$, $t>0$,  is the inverse of $\mathpzc{H}^\nu (t)$, $t>0$ and possesses distribution $\mathpzc{l}_\nu (x, t)$, $x \geq 0$, $t>0$.
	\item[$\bullet$] $T(t)$, $t>0$, is a telegraph process with parameters $c>0$ and $\lambda >0$ and law $p_T (x, t)$, $-ct<x<ct$, $t>0$.
	\item[$\bullet$]	$\bm{W}_n \l t \r \, = \, \bm{S}_n^{2\beta} \l c^2 \mathpzc{L}^\nu (t) \r$, $t>0$, has law $w_\nu^\beta \l \bm{x}, t \r$, $\bm{x} \in \mathbb{R}^n$, $t>0$.
	\item[$\bullet$]   $\mathpzc{W} (t) = T \l | B (t) | \r$, $t>0$, has distribution $\mathpzc{w} (x, t)$, $x \in \mathbb{R}$, $t>0$.
	\item[$\bullet$]   $\bm{T} (t)$, $t>0$, is the planar process with infinite directions, parameters $c, \lambda >0$ and law $r(x, y, t)$, $(x,y) \in C_{ct} = \ll (x, y) \in \mathbb{R}^2 : x^2+y^2 < c^2t^2 \rr$, $t>0$.
	\item[$\bullet$]   $\bm{\mathfrak{T}} (t)$, $t>0$, is the planar process with infinite directions, parameters $c, \lambda >0$ and law $\mathfrak{r}(x, y, t)$, $(x,y) \in C_{ct} = \ll (x, y) \in \mathbb{R}^2 : x^2+y^2 < c^2t^2 \rr$, $t>0$, constructed by disregading displacements started off only by even-labelled Poisson events.
\item[$\bullet$] $\bm{Q}(t) = \bm{T} \l \left| B(t) \right| \r$, $t>0$, has law $q(x, y, t)$, $(x, y) \in \mathbb{R}^2$, $t>0$.

	\item[$\bullet$] By $\widetilde{f}$ we denote the Laplace transform of the function $f$ and by $\widehat{f}$ we denote its Fourier transform.
	\end{enumerate}

	\subsection{Preliminaries}
	\label{prelimin}
	Let us consider a stable process $S^\nu (t)$, $t>0$, $0<\nu\leq 2$, $\nu \neq 1$, with characteristic function
	\begin{equation} 
		\mathbb{E}e^{i \xi S^{\nu}(t)} \, = \, e^{-\sigma \left| \xi \right|^\nu t \l 1-i\theta \, \textrm{sign}(\xi) \, \tan \frac{\nu \pi}{2} \r}
		\label{stable caratt}
	\end{equation}
where $\theta \in \left[ -1,1 \right]$ is the skewness parameter and
\begin{equation} 
	\sigma \, = \, \cos \frac{\pi \nu}{2}.
	\label{}
\end{equation}
 For $\theta = 1$ the distribution corresponding to \eqref{stable caratt} is totally positively skewed and for $\theta = -1$ is totally negatively skewed.
The stable process with stationary and independent increments, totally positively skewed will be denoted as $H^\nu (t)$, $t>0$. We note that the density $h_\nu (x, t)$, of $H^{\nu}(t)$, is zero at $x=0$ as the following calculation show
\begin{align}
	h_\nu (0, t) \, & = \, \frac{1}{2\pi} \int_{-\infty}^\infty \mathbb{E}e^{i\xi H^\nu (t)} \, d\xi \, = \, \frac{1}{2\pi} \int_{-\infty}^\infty e^{-\sigma \left| \xi \right|^\nu t \l 1-i\tan \frac{\nu \pi }{2} \r} d\xi \notag \\
	& = \, \frac{1}{2\pi} \left[ \int_0^\infty e^{-\sigma \left| \xi \right|^\nu t \l 1-i\tan \frac{\nu \pi }{2} \r} d\xi + \int_{-\infty}^0 e^{-\sigma \left| \xi \right|^\nu t \l 1+i\tan \frac{\nu \pi }{2} \r} d\xi  \right] \notag \\
	& = \, \frac{1}{2\pi} \left[ \int_0^\infty e^{-\left| \xi \right|^\nu t e^{ -\frac{i\nu \pi}{2} }} d\xi + \int_0^\infty e^{-\left| \xi \right|^\nu t e^{ \frac{i \nu \pi}{2} }} d\xi \right] \notag \\
	& = \, \frac{1}{2\pi} \left[ \int_0^\infty e^{-z} \l \frac{z}{t} \r^{ \frac{1}{\nu}-1 } e^{ \frac{i\pi}{2} } dz + \int_0^\infty e^{-z} \l \frac{z}{t} \r^{ \frac{1}{\nu}-1 } \frac{1}{t} e^{-\frac{i\pi}{2}  } dz \right] \notag \\
	& = \, \frac{\cos \frac{\pi}{2}}{\pi} \int_0^\infty e^{-z} \l \frac{z}{t} \r^{ \frac{1}{\nu}-1 } \frac{1}{t}  \, dz \, = \, 0.
	\label{}
\end{align}
The positively skewed stable r.v. $H^\nu (t)$ has $x$-Laplace transform
\begin{equation} 
	\widetilde{h_\nu} \l \mu, t \r \, = \,  \mathbb{E}e^{-\mu H^\nu(t)} \, = \, e^{-t\mu^\nu}, \qquad 0<\nu <1,
	\label{laplace di h}
\end{equation}
and therefore Fourier transform
\begin{align}
	\widehat{h_\nu} \l \xi, t \r \, = \, \mathbb{E}e^{i\xi H^\nu (t)} \, & = \, \mathbb{E} \l e^{- \l -i \xi \r H^\nu (t)} \r \, = \, e^{-t \l \left| \xi \right| e^{- \frac{i \pi}{2} \textrm{sign}(\xi)} \r^\nu} \notag \\
	& =  \, e^{-t \left| \xi \right|^\nu \cos \frac{\pi \nu}{2}  \l 1-i \, \textrm{sign} (\xi) \tan \frac{\pi \nu}{2} \r }.
	\label{fourier subordinatore}
\end{align}
This shows once again that the skeweness parameter is $\theta =1$. 

The probability law $h_\nu (x, t)$, of $H^\nu (t)$, $t>0$, solves the problem
\begin{align}
	\begin{cases}
	\l \frac{\partial}{\partial t} + \frac{\partial^\nu}{\partial x^\nu} \r h_\nu (x, t) \, = \, 0, \qquad & x >0, t>0, 0<\nu<1, \\
	h_\nu (0, t) \, = \, 0, \\
	h_\nu (x, 0) \, = \, \delta (x).
	\end{cases}
	\label{bifp}
\end{align}
By taking the $x$-Laplace transform of the Riemann-Liouville fractional derivative appearing in \eqref{bifp} we have that
\begin{align}
	& \mathcal{L}  \left[ \frac{\partial^\nu}{\partial x^\nu} h_\nu (x, t) \right] (\mu) \,  =  \,  \int_0^\infty e^{-\mu x} \frac{\partial^\nu}{\partial x^\nu} h_\nu (x, t) \, dx \notag \\
	 = \, &\int_0^\infty e^{-\mu x} \left[ \frac{1}{\Gamma \l 1-\nu \r} \frac{d}{dx} \int_0^x \frac{h_\nu (z, t)}{\l x-z \r^\nu} dz  \right]dx \notag \\
	 = \, & \int_0^\infty e^{-\mu x} \left[ \frac{1}{\Gamma \l 1-\nu \r} \int_0^x \frac{d}{dx} \frac{  h_\nu (x-z, t)}{z^\nu} dz + \frac{h_\nu (0, t)}{\Gamma \l 1-\nu \r \, x^\nu} \right] dx \notag \\
	 = \, &\frac{h_\nu (0, t)}{\Gamma \l 1-\nu \r} \int_0^\infty e^{-\mu x} x^{1-\nu-1} \, dx   + \frac{1}{\Gamma \l 1-\nu \r} \int_0^\infty \frac{dz}{z^\nu} \int_z^\infty dx \, e^{-\mu x} \frac{d}{dx} h_\nu (x-z, t)  \notag \\
	 = \, & h_\nu (0, t) \mu^{\nu-1} + \frac{1}{\Gamma \l 1-\nu \r} \int_0^\infty e^{-\mu z} z^{-\nu} dz \int_0^\infty e^{-\mu x} \frac{d}{dx} h_\nu (x, t) dx \notag \\
	 = \, & h_\nu (0, t) \mu^{\nu -1} + \left[ \int_0^\infty e^{-\mu x} h_\nu (x, t) dx \right] \mu \, \frac{1}{\mu^{1-\nu}} - \mu^{\nu-1} h_\nu (0, t) \,
	= \,  \mu^\nu  \widetilde{h_\nu} \l \mu, t \r.
	\label{laplace riemann}
\end{align}
Therefore
\begin{align}
	\begin{cases}
	\frac{\partial}{\partial t} \widetilde{h_\nu} \l \mu, t \r + \mu^\nu \widetilde{h_\nu} \l \mu, t \r \, = \, 0, \qquad \mu >0, t>0, \\
	\widetilde{h_\nu} \l \mu, 0 \r \, = \, 1,
	\end{cases}
	\label{}
\end{align}
so that
\begin{equation} 
	\widetilde{h_\nu} \l \mu, t \r \, = \, e^{-\mu^\nu t}.
	\label{153}
\end{equation}
In other words the density of a positively skewed stable r.v. solves the space-fractional problem \eqref{bifp}.

We will also deal with the inverse process of $H^\nu (t)$, $t>0$, say $L^\nu (t)$, $t>0$, for which
\begin{equation} 
	\Pr \ll H^\nu (x) > t \rr \, = \, \Pr \ll L^\nu (t) < x \rr, \qquad x>0, t>0.
	\label{}
\end{equation}
Such a process has non-negative, non-stationary and non-independent increments. Furthemore we recall that the law $l_\nu (x, t)$ of $L^\nu (t)$, can be written as
\begin{equation} 
	l_\nu (x, t) \, = \, \frac{1}{t^\nu} W_{-\nu, 1-\nu} \l -\frac{x}{t^\nu} \r, \qquad x \geq 0, t>0,
	\label{155}
\end{equation}
where 
\begin{equation} 
	W_{a, b } (x) \, = \, \sum_{k=0}^\infty \frac{x^k}{k! \, \Gamma \l ak+b \r}, \qquad x \in \mathbb{R}, a> -1, b \in \mathbb{C},
\end{equation}
is the Wright function, and has Laplace transform
\begin{align}
		\widetilde{l_\nu} (x, \mu) \, = \, \int_0^\infty \, e^{-\mu t} l_\nu (x, t) dt \, = \, \int_0^\infty \, e^{-\mu t}\frac{1}{t^\nu} W_{-\nu, 1-\nu} \l -\frac{x}{t^\nu} \r dt \, = \, \mu^{\nu -1} e^{-x \mu^\nu}.
		\label{lapl l}
		\end{align}
%

\section{Sum of stable subordinators, $\mathpzc{H}^\nu (t) = H_1^{2\nu}(t) + (2\lambda)^{\frac{1}{\nu}} H_2^\nu (t)$}

For the construction of the vector process $\bm{W}_n(t)= \bm{S}_n^{2\beta} \l c^2 \mathpzc{L}^\nu (t) \r$, $t>0$, whose distribution is driven by the general space-time fractional telegraph equation \eqref{1.1}, we need the sum $\mathpzc{H}^\nu (t)$, $t>0$, of two independent positively skewed processes. The second step consists in constructing the process $\mathpzc{L}^\nu (t)$, $t>0$, inverse to $\mathpzc{H}^\nu (t)$, $t>0$. We now start by considering the following sum
\begin{equation} 
	\mathpzc{H}^\nu(t) \, = \, H_1^{2\nu} (t) + \l 2\lambda \r^{ \frac{1}{\nu} } H_2^{\nu} (t), \qquad t>0, 0<\nu \leq \frac{1}{2},
	\label{}
\end{equation}
with $H_1^{2\nu}$, $H_2^{\nu}$, independent, positively-skewed, stable random variables, $\lambda >0$.
The distribution of $\mathpzc{H}^\nu (t)$ can be written as
\begin{align}
	\mathpzc{h}_\nu \l x, t \r \, = \, \int_0^x h_{2\nu} (y, t) \, h_{\nu} (x-y, 2\lambda t) \, dy.
	\label{convoluzione}
\end{align}
Taking the double Laplace transform of \eqref{convoluzione}, with respect to $t$ and $x$, we get
\begin{align}
	\widetilde{\widetilde{\mathpzc{h}_\nu}} \l \gamma, \mu \r  \, & = \, \int_0^\infty e^{-\mu t}  \int_0^\infty e^{-\gamma x} \mathpzc{h}_\nu (x, t) \, dx \, dt \, = \, \int_0^\infty e^{-\mu t-t\gamma^{2\nu}-2\lambda t \gamma^\nu} dt \notag \\
	& = \, \frac{1}{\gamma^{2\nu} + 2\lambda \gamma^\nu + \mu} \, = \, \left[ \frac{1}{\gamma^\nu - r_2} - \frac{1}{\gamma^\nu - r_1} \right] \frac{1}{r_2-r_1}
	\label{xtlapl}
\end{align}
where, for $0 < \mu < \lambda^2$,
\begin{align}
	\begin{cases}
	r_1 \, = \, -\lambda - \sqrt{\lambda^2 - \mu}, \\
	r_2 \, = \, -\lambda + \sqrt{\lambda^2 -\mu}.
	\end{cases}
	\label{}
\end{align}
By means of formula
\begin{equation} 
	\int_0^\infty e^{-\gamma x } \, x^{\alpha -1} \, E_{\alpha, \alpha} \l \eta x^\alpha \r \, dx \, = \, \frac{1}{\gamma^\alpha -\eta},
	\label{}
\end{equation}
where $E_{\nu, \nu} (z)$ is the Mittag-Leffler function defined in \eqref{def mittag leffler}, we can invert the $x$-Laplace transform in \eqref{xtlapl} obtaining, for $\mu < \lambda^2$,
\begin{align}
	& \widetilde{\mathpzc{h}_\nu} \l x, \mu \r \,   = \notag \\
	 = \, & \frac{x^{\nu -1}}{2\sqrt{\lambda^2 - \mu}} \left[ E_{\nu , \nu} \l \l -\lambda  + \sqrt{\lambda^2 - \mu} \r x^\nu \r  - E_{\nu , \nu} \l \l -\lambda - \sqrt{\lambda^2 - \mu} \r x^\nu \r \right] \notag \\
	  = \, & \frac{1}{2\sqrt{\lambda^2 - \mu}} \left[ \frac{1}{-\lambda + \sqrt{\lambda^2 - \mu}} \frac{\partial}{\partial x} E_{\nu, 1} \l \l -\lambda + \sqrt{\lambda^2 - \mu} \r x^\nu \r  \,  \right. \notag \\
	&  \left. - \frac{1}{-\lambda - \sqrt{\lambda^2 - \mu}} \frac{\partial}{\partial x} E_{\nu, 1} \l \l -\lambda - \sqrt{\lambda^2 - \mu} \r x^\nu \r \right].
	\label{above}
\end{align}
Formula \eqref{above} gives the explicit form of the $t$-Laplace transform of $\mathpzc{h}_\nu (x, t)$ in terms of Mittag-Leffler functions.
In view of formula
\begin{equation} 
	E_{\nu, 1} \l -\lambda t^\nu \r \, = \, \frac{1}{\pi} \int_0^\infty \frac{e^{-\lambda^{ \frac{1}{\nu} }tx}x^{\nu-1} \sin \pi \nu}{x^{2\nu}+1+2x^\nu \cos \pi \nu} \, dx, \qquad 0 < \nu < 1,
	\label{}
\end{equation}
we have that
\begin{align}
	\widetilde{\mathpzc{h}_\nu} \l x, \mu \r \,  = \, & \frac{1}{2\sqrt{\lambda^2 - \mu}} \left[  \frac{1}{-\lambda + \sqrt{\lambda^2 - \mu}} \frac{\partial}{\partial x}   \int_0^\infty \frac{e^{-xy \l \lambda - \sqrt{\lambda^2 - \mu} \r^{ \frac{1}{\nu} }}y^{\nu-1} \sin \pi \nu \, dy}{\pi \l y^{2\nu}+1+2y^\nu \cos \pi \nu \r } \right.   \notag \\
	&  + \left.  \frac{1}{\lambda + \sqrt{\lambda^2 - \mu}} \frac{\partial}{\partial x}  \int_0^\infty \frac{e^{-xy \l \lambda  + \sqrt{\lambda^2 - \mu} \r^{ \frac{1}{\nu} }}y^{\nu -1} \sin \pi \nu \, dy}{\pi \l y^{2\nu}+1+2y^\nu \cos \pi \nu \r} \right] \notag \\
	 = \, & \int_0^\infty \frac{dy \; y^\nu \sin \pi \nu}{\pi \l y^{2\nu} + 1 + 2y^\nu \cos \pi \nu \r} \frac{1}{2\sqrt{\lambda^2 - \mu}} \Bigg[ \l \lambda - \sqrt{\lambda^2 - \mu} \r^{ \frac{1}{\nu}-1 }  \, \bm{\cdot} \notag \\
	&  \bm{\cdot}  e^{-xy \l \lambda - \sqrt{\lambda^2 - \mu} \r^{ \frac{1}{\nu} }} - \l \lambda + \sqrt{\lambda^2 - \mu} \r^{ \frac{1}{\nu}-1 } e^{-xy \l \lambda + \sqrt{\lambda^2 - \mu} \r^{ \frac{1}{\nu} }} \Bigg] \notag \\
	  = \, & \mathbb{E} \ll \frac{\mathpzc{U}^\nu}{2\sqrt{\lambda^2 - \mu}} \left[ \l -r_2 \r^{\frac{1}{\nu}-1} e^{-x \mathpzc{U}^\nu (-r_2)^{\frac{1}{\nu}}} - \l -r_1 \r^{\frac{1}{\nu}-1} e^{-x \mathpzc{U}^\nu (-r_1)^{\frac{1}{\nu}}} \right] \rr \notag \\
	 = \, & \frac{1}{r_2-r_1} \frac{\partial}{\partial x} \mathbb{E} \left[ \frac{e^{-x \mathpzc{U}^\nu (-r_2)^{\frac{1}{\nu}}}}{r_2} - \frac{e^{-x \mathpzc{U}^\nu (-r_1)^{\frac{1}{\nu}}}}{r_1} \right],
	\label{}
\end{align}
where $\mathpzc{U}^\nu$ is the Lamperti distribution with density
\begin{equation} 
	 \frac{\Pr \ll \mathpzc{U^\nu} \in d\mathpzc{u} \rr}{d\mathpzc{u}} \, = \, \frac{\sin \pi \nu}{\pi} \frac{\mathpzc{u}^{\nu -1}}{1+ \mathpzc{u}^{2\nu} + 2 \mathpzc{u}^\nu \cos \pi \nu}, \qquad \mathpzc{u} >0,
	\label{}
\end{equation}
and represents the law of the ratio of two independent stable r.v.'s of the same order $\nu$.
\begin{te}
The law $\mathpzc{h}_\nu (x, t)$ of the process $\mathpzc{H}^\nu (t) \, = \, H_1^{2\nu} (t) + \l 2\lambda \r^{ \frac{1}{\nu} } H_2^\nu (t)$ solves the fractional problem
\begin{align}
	\begin{cases}
	\frac{\partial}{\partial t} \mathpzc{h}_\nu \l x, t \r \, = \, - \l \frac{\partial^{2\nu}}{\partial x^{2\nu}} +  2\lambda  \frac{\partial^\nu}{\partial x^\nu} \r \mathpzc{h}_\nu (x, t), \qquad x >0, t>0, 0 < \nu < \frac{1}{2}, \\
	\mathpzc{h}_\nu (0, t) \, = \, 0, \\
	\mathpzc{h}_\nu (x, 0) \, = \, \delta (x).
	\label{governing equation}
	\end{cases}
\end{align}
The fractional derivatives appearing in \eqref{governing equation} are intended in the Riemann-Liouville sense.
\begin{proof}
By considering \eqref{fourier subordinatore}, we have that the Fourier transform of $\mathpzc{h}_\nu (x,t)$ is written as
\begin{align}
	\widehat{\mathpzc{h}_\nu} (\xi, t) \, & = \, \mathbb{E}e^{i\xi \mathpzc{H}^\nu (t)}  \, = \, \mathbb{E}e^{i\xi \left[ H^{2\nu}(t) + \l 2\lambda \r^{ \frac{1}{\nu} } H^\nu (t) \right]}  \, = \, \mathbb{E}e^{i\xi H^{2\nu}(t)} e^{i\xi H^\nu (2\lambda t)}  \notag \\
	& = \, e^{-t \left| \xi \right|^{2\nu} \cos \pi \nu \l 1-i \, \textrm{sign}(\xi) \tan \pi \nu \r - 2\lambda t\left| \xi \right|^\nu \cos \frac{\pi \nu}{2} \l 1-i \, \textrm{sign}(\xi) \tan \frac{\pi \nu}{2} \r} \notag \\
	& = \, e^{-t \l \left| \xi \right| e^{- \frac{i\pi}{2} \, \textrm{sign}(\xi)} \r^{2\nu}-2\lambda t \l \left| \xi \right| e^{ -\frac{i\pi}{2} \textrm{sign}(\xi) }  \r^\nu},
	\label{}
\end{align}
and thus
\begin{align}
	\frac{\partial}{\partial t} \widehat{\mathpzc{h}_\nu} \l \xi, t \r \,  = \, & \left[ - \l \left| \xi \right| e^{- \frac{i \pi}{2} \, \textrm{sign}(\xi) } \r^{2\nu} - 2\lambda \l \left| \xi \right| e^{-\frac{i \pi}{2} \, \textrm{sign}(\xi) } \r^\nu \right] \, \bm{\cdot} \notag \\
	&  \bm{\cdot} \,e^{-t \l \left| \xi \right| e^{- \frac{i\pi}{2} \, \textrm{sign}(\xi)} \r^{2\nu}-2\lambda t \l \left| \xi \right| e^{ -\frac{i\pi}{2} \textrm{sign}(\xi) }  \r^\nu}
	\label{derivative of fourier}
\end{align}
In view of the relationship
\begin{equation} 
	\left| \xi \right| e^{-\frac{i \pi}{2} \, \textrm{sign}(\xi)} \, = \, -i \xi
	\label{}
\end{equation}
we have that formula \eqref{derivative of fourier} can be rewritten as
\begin{align}
	\frac{\partial}{\partial t} \widehat{\mathpzc{h}_\nu} \l \xi, t \r \, & = \, \left[ - \l -i\xi \r^{2\nu} - 2\lambda \l - i \xi \r^\nu \right]   e^{-t \l -i \xi \r^{2\nu} - 2\lambda t \l -i \xi \r^\nu}.
	\label{derivative final}
\end{align}
In \eqref{laplace riemann} we have shown that
\begin{align}
	\mathcal{L} \left[ \frac{\partial^\nu}{\partial x^\nu} h_\nu \l x, t \r \right] (\mu) \, = \, \int_0^\infty e^{-\mu x} \frac{\partial^\nu}{\partial x^\nu} h_\nu (x, t) \, dx \, = \,  \mu^\nu \widetilde{h_\nu} (\mu, t)
	\label{}
\end{align}
and thus for a sufficiently good function $f$ we have the following Fourier transform
\begin{align}
	\mathcal{F} \left[ \frac{\partial^\nu}{\partial x^\nu} f(x) \right] (\xi) \, & = \, \int_0^\infty e^{- \l -i \xi \r x} \frac{\partial^\nu}{\partial x^{\nu}} f(x) \, dx \, = \,  \l -i\xi \r^\nu \widehat{f}(\xi).
	\label{fourier riemann liouville}
\end{align}
In view of \eqref{fourier riemann liouville} we have that the Fourier transform of the right-hand side of the equation \eqref{governing equation}, equipped with the boundary conditions, is written as
\begin{align}
	& -\mathcal{F} \left[ \frac{\partial^{2\nu}}{\partial x^{2\nu}} \mathpzc{h}_\nu \l x, t \r +  2\lambda  \frac{\partial^\nu}{\partial x^\nu} \mathpzc{h}_\nu \l x, t \r \right] \l \xi \r \,  = \notag \\
	 = \, & -\int_0^\infty e^{-\l-i \xi \r x }  \frac{\partial^{2\nu}}{\partial x^{2\nu}} \mathpzc{h}_\nu \l x, t \r \, dx \, -   2\lambda  \int_0^\infty e^{- \l -i \xi \r x} \frac{\partial^\nu}{\partial x^\nu} \mathpzc{h}_\nu \l x, t \r \, dx \notag \\
	 = \, &- \l  \l -i\xi \r^{2\nu}  + 2\lambda  \l -i\xi \r^{\nu} \r \widehat{\mathpzc{h}_\nu} \l \xi, t \r \notag \\
	 = \, &-  \l  \l -i\xi \r^{2\nu}  + 2\lambda  \l -i\xi \r^{\nu} \r  \, e^{-t \l \left| \xi \right| e^{- \frac{i\pi}{2} \, \textrm{sign}(\xi)} \r^{2\nu}-2\lambda t \l \left| \xi \right| e^{ -\frac{i\pi}{2} \textrm{sign}(\xi) }  \r^\nu} \notag \\
	 = \, & -  \l  \l -i\xi \r^{2\nu}  + 2\lambda \l -i\xi \r^{\nu} \r   \, e^{- t \l -i \xi \r^{2\nu} - 2\lambda t \l -i \xi \r^\nu},
	\label{}
\end{align}
which coincides with formula \eqref{derivative final}. This is tantamount to saying that the Fourier transform $\widehat{\mathpzc{h}_\nu} \l \xi, t \r$ is the solution to
\begin{align}
	\begin{cases}
	\frac{\partial}{\partial t} \widehat{\mathpzc{h}_\nu} \l \xi, t \r \, = \, - \l  \l -i\xi \r^{2\nu}  + 2\lambda  \l -i\xi \r^{\nu} \r \widehat{\mathpzc{h}_\nu} \l \xi, t \r, \qquad \xi \in \mathbb{R}, t>0, \\
	\widehat{\mathpzc{h}_\nu} \l \xi, 0 \r \, = \, 1,
	\end{cases}
	\label{}
\end{align}
and this completes the proof.
\end{proof}
\end{te}

\subsection{The inverse process $\mathpzc{L}^\nu (t)$}
Let $\mathpzc{L}^\nu (t)$, $t>0$, be the inverse process of $\mathpzc{H}^\nu (t)$, $t>0$, as defined in \eqref{119} for which
\begin{equation} 
	\Pr \ll \mathpzc{L}^\nu (t) < x \rr \, = \, \Pr \ll \mathpzc{H}^\nu (x) > t \rr, \qquad x, t >0,
	\label{construction L}
\end{equation}
and let $\mathpzc{l}_\nu (x, t)$ be the law of $\mathpzc{L}^\nu (t)$, $t>0$. We have the following result.
\begin{te}
\label{teorema inverso}
The law $\mathpzc{l}_\nu(x, t)$ of the process $\mathpzc{L}^\nu(t)$, $t>0$, solves the time-fractional boundary-initial problem
\begin{align}
	\begin{cases}
	\l \frac{\partial^{2\nu}}{\partial t^{2\nu}} + 2\lambda \frac{\partial^\nu}{\partial t^\nu} \r \mathpzc{l}_\nu (x, t) \, = \, -\frac{\partial}{\partial x} \mathpzc{l}_\nu (x, t), \qquad x>0, t>0, 0<\nu < \frac{1}{2}, \\
	\mathpzc{l}_\nu (x, 0) \, = \, \delta (x), \\
	\mathpzc{l}_\nu (0, t) \, = \, \frac{t^{-2\nu}}{\Gamma \l 1-2\nu \r} + 2\lambda \frac{t^{-\nu}}{\Gamma \l 1-\nu \r},
	\end{cases}
	\label{inverse problem}
\end{align}
and has $x$-Laplace transform which reads, for $0<\gamma<\lambda^2$,
\begin{align}
	\widetilde{\mathpzc{l}_\nu} \l \gamma, t \r \, = \, \frac{1}{2} \left[ \l 1+ \frac{\lambda}{\sqrt{\lambda^2- \gamma}} \r E_{\nu,1} \l r_1 t^\nu \r + \l 1-\frac{\lambda}{\sqrt{\lambda^2 - \gamma}} \r E_{\nu,1} \l r_2 t^\nu \r  \right],
	\label{laplace di lbold}
\end{align}
where 
\begin{equation} 
	r_1 \, = \, -\lambda + \sqrt{\lambda^2 - \gamma}, \qquad r_2 \, = \, -\lambda -\sqrt{\lambda^2 - \gamma}.
	\label{}
\end{equation}
The fractional derivatives appearing in \eqref{inverse problem} are intended in the Riemann-Liouville sense.
\begin{proof}
We first show that the analytical solution to the problem \eqref{inverse problem} has double Laplace transform $\widetilde{\widetilde{\mathpzc{l}_\nu}}\l \gamma, \mu \r$ written as
\begin{equation} 
	\widetilde{\widetilde{\mathpzc{l}_\nu}}\l \gamma, \mu \r \, = \, \frac{\mu^{2\nu-1}+2\lambda \mu^{\nu-1}}{\mu^{2\nu}+2\lambda \mu^\nu + \gamma}.
	\label{fourier laplace analytic}
\end{equation}
By taking the $t$-Laplace transform of the equation in \eqref{inverse problem} we have that
\begin{align}
	\mu^{2\nu} \widetilde{\mathpzc{l}_\nu} \l x, \mu \r + 2\lambda \mu^\nu \widetilde{\mathpzc{l}_\nu} \l x, \mu \r \, = \, -\frac{\partial}{\partial x} \widetilde{\mathpzc{l}_\nu} \l x, \mu \r.
	\label{tlaplace inverse problem}
\end{align}
By taking into account the boundary condition and performing the $x$-Laplace transform of \eqref{tlaplace inverse problem} we have that
\begin{equation} 
	\l \mu^{2\nu} + 2\lambda \mu^\nu \r \widetilde{\widetilde{\mathpzc{l}_\nu}} \l \gamma, \mu \r \, = \, \widetilde{\mathpzc{l}_\nu} \l 0, \mu \r - \gamma \, \widetilde{\widetilde{\mathpzc{l}_\nu}} \l \gamma, \mu \r.
	\label{}
\end{equation}
Now, by considering the boundary condition, we get that
\begin{align}
	\widetilde{\mathpzc{l}_\nu} \l 0, \mu \r \, & = \, \int_0^\infty dt \, e^{-\mu t} \mathpzc{l}_\nu \l 0, t \r \, = \, \int_0^\infty dt \, e^{-\mu t} \left[ \frac{t^{-2\nu}}{\Gamma \l 1-2\nu \r} + 2\lambda \frac{t^{-\nu}}{\Gamma \l 1-\nu \r} \right] \notag \\
	& = \, \mu^{2\nu -1} + 2\lambda \mu^{\nu -1},
	\label{}
\end{align}
and thus
\begin{align}
	 \widetilde{\widetilde{\mathpzc{l}_\nu}} \l \gamma, \mu \r \, = \, \frac{\mu^{2\nu -1} + 2\lambda \mu^{\nu -1}}{\mu^{2\nu}+2\lambda \mu^{\nu} + \gamma }.
	\label{}
\end{align}
Now we show that the double Laplace transform of the law $\mathpzc{l}_\nu \l x, t \r$ coincides with \eqref{fourier laplace analytic}.
We first recall that
\begin{align}
	\widetilde{\mathpzc{h}_\nu} (\mu, x) \, & = \, \int_0^\infty dt \, e^{-\mu t} \mathpzc{h}_\nu (t, x) \, = \,  \mathbb{E} e^{-\mu \mathpzc{H}^\nu (x)} \, = \, \mathbb{E} e^{-\mu H^{2\nu}(x)} \mathbb{E}e^{-\mu H^\nu \l 2\lambda x\r} \notag \\
	& = \, \widetilde{h_{2\nu}} \l \mu, x \r \, \widetilde{h_{\nu}} \l \mu, 2\lambda x \r \, = \, e^{-x\mu^{2\nu} - x 2\lambda \mu^{\nu}}, \qquad x >0,
	\label{lapl proc}
\end{align}
where we used result \eqref{laplace di h}.
 By considering the construction of the process $\mathpzc{L}^\nu (t)$, $t>0$, as the inverse process of $\mathpzc{H}^\nu (t)$, $t>0$, as stated in \eqref{construction L}, we get
\begin{align}
	\mathpzc{l}_\nu (x, t)  \, = \, \frac{\Pr \ll \mathpzc{L}^\nu (t) \in dx \rr}{dx}  \, = \, -\frac{\partial}{\partial x} \Pr \ll \mathpzc{H}^\nu (x) < t \rr \, = \, - \frac{\partial}{\partial x} \int_0^t \mathpzc{h}_\nu (s, x) \, ds.
	\label{density of l}
\end{align}
In view of \eqref{density of l}, the double Laplace transform of $\mathpzc{l}_\nu (x, t)$ can be obtained observing that
\begin{align}
	\widetilde{\widetilde{\mathpzc{l}_\nu}} \l \gamma, \mu \r \, & = \, \int_0^\infty dx \, e^{-\gamma x} \int_0^\infty dt \, e^{-\mu t} \left[- \frac{\partial}{\partial x} \int_0^t \mathpzc{h}_\nu (s, x) \, ds \right]   \notag \\
	& = \,- \int_0^\infty dx \, e^{-\gamma x} \frac{\partial}{\partial x} \int_0^\infty dt \, e^{-\mu t} \int_0^t \mathpzc{h}_\nu (s, x) \, ds \notag \\
		& = \, -\frac{1}{\mu} \int_0^\infty dx \, e^{-\gamma x}  \frac{\partial}{\partial x} \widetilde{\mathpzc{h}_\nu} \l x, \mu \r  \, = \,  -\frac{1}{\mu} \int_0^\infty dx \, e^{-\gamma x} \left[ \frac{\partial}{\partial x} e^{-x\mu^{2\nu}-2\lambda x \mu^\nu } \right] \notag \\
	& = \, \l \mu^{2\nu-1} + 2\lambda \mu^{\nu-1} \r \int_0^\infty dx \, e^{-\gamma x-x\mu^{2\nu}-2\lambda x \mu^\nu} \, = \, \frac{\mu^{2\nu-1 }+2\lambda \mu^{\nu-1}}{\mu^{2\nu}+ 2\lambda \mu^\nu  + \gamma },
	\label{}
\end{align}
which coincides with \eqref{fourier laplace analytic}.
Now we pass to the derivation of the $x$-Laplace transform of $\mathpzc{l}_\nu \l x, t \r$. We can write
\begin{align}
	\widetilde{\widetilde{\mathpzc{l}_\nu}} \l \gamma, \mu \r \, & = \, \frac{\mu^{2\nu-1 }+2\lambda \mu^{\nu-1}}{\mu^{2\nu}+ 2\lambda \mu^\nu  + \gamma } \, = \, \frac{\mu^{\nu-1}}{\mu^\nu - r_1} + \frac{\mu^{\nu-1}}{\mu^\nu - r_2} - \frac{\mu^{2\nu-1}}{\l \mu^\nu - r_1 \r \l \mu^\nu - r_2 \r} \notag \\
	& = \, \frac{\mu^{\nu-1}}{\mu^\nu - r_1} + \frac{\mu^{\nu-1}}{\mu^\nu - r_2} - \left[ \frac{\mu^{\nu- (1-\nu)}}{\mu^\nu - r_1} - \frac{\mu^{\nu-(1-\nu)}}{\mu^\nu - r_2} \right] \, \frac{1}{2\sqrt{\lambda^2 - \gamma}},
	\label{double laplace di lnu}
\end{align}
where
\begin{equation} 
	r_1 \, = \, -\lambda + \sqrt{\lambda^2 - \gamma}, \qquad r_2 \, = \, -\lambda -\sqrt{\lambda^2 - \gamma}.
	\label{}
\end{equation}
Now we need the following results
\begin{align}
	&\int_0^\infty e^{-\mu t} E_{\nu,1} \l r_j t^\nu \r \, dt \, = \, \frac{\mu^{\nu-1}}{\mu^\nu - r_j}, \qquad j = 1, 2, \notag \\
	& \int_0^\infty e^{-\mu t} t^{(1-\nu)-1} E_{\nu, 1-\nu} \l r_j t^\nu \r \, dt \, = \, \frac{\mu^{2\nu-1}}{\mu^\nu - r_j}.
	\label{proprietàinversione}
\end{align}
Therefore
\begin{align}
	\widetilde{\mathpzc{l}_\nu} \l \gamma, t \r \, & = \, E_{\nu,1} \l r_1 t^\nu \r + E_{\nu, 1} \l r_2 t^\nu \r - \frac{t^{-\nu}}{2\sqrt{\lambda^2 - \gamma}} \left[  E_{\nu, 1-\nu} \l  r_1 t^\nu  \r -  E_{\nu, 1-\nu} \l r_2 t^\nu \r \right].
	\label{}
\end{align}
Since
\begin{align}
	E_{\nu, 1-\nu} \l z \r \, = \, z E_{\nu,1} (z) + \frac{1}{\Gamma \l 1-\nu \r}
	\label{}
\end{align}
we have that
\begin{align}
	\widetilde{\mathpzc{l}_\nu} \l \gamma, t \r \, & = \,  E_{\nu,1} \l r_1 t^\nu \r + E_{\nu, 1} \l r_2 t^\nu \r - \frac{t^{-\nu}}{2\sqrt{\lambda^2 - \gamma}} \left[ r_1 t^\nu E_{\nu,1} \l r_1 t^\nu \r - r_2 t^\nu E_{\nu,1} \l r_2 t^\nu \r \right] \notag \\
	& = \, \l 1- \frac{-\lambda + \sqrt{\lambda^2 - \gamma}}{2\sqrt{\lambda^2 - \gamma}} \r E_{\nu,1} \l r_1 t^\nu \r + \l 1- \frac{\lambda + \sqrt{\lambda^2 - \gamma}}{2\sqrt{\lambda^2 - \gamma}} \r E_{\nu , 1} \l r_2 t^\nu \r \notag \\
	& = \, \frac{1}{2} \left[ \l 1+ \frac{\lambda}{\sqrt{\lambda^2 - \gamma}} \r E_{\nu, 1} \l r_1 t^\nu \r + \l 1- \frac{\lambda}{\sqrt{\lambda^2 - \gamma}} \r E_{\nu, 1} \l r_2 t^\nu \r   \right],
	\label{236}
\end{align}
which coincides with \eqref{laplace di lbold}.

Now we check that the Laplace transform \eqref{236} solves the fractional equation
\begin{align} 
	\l \frac{\partial^{2\nu}}{\partial t^{2\nu}}  +2\lambda \frac{\partial^\nu}{\partial t^\nu} \r 	\widetilde{\mathpzc{l}_\nu} \l \gamma, t \r \, & = \, -\gamma \, 	\widetilde{\mathpzc{l}_\nu} \l \gamma, t \r + \mathpzc{l}_\nu \l 0, t \r \notag \\
	& = \,  -\gamma \, 	\widetilde{\mathpzc{l}_\nu} \l \gamma, t \r + \frac{t^{-2\nu}}{\Gamma (1-2\nu)} + 2\lambda \frac{t^{-\nu}}{\Gamma \l 1- \nu \r}
	\label{stellina}
\end{align}
which is the $x$-Laplace transform of the equation appearing in \eqref{inverse problem}.
Since
\begin{align}
	&\frac{\partial^{2\nu}}{\partial t^{2\nu}} \widetilde{\mathpzc{l}_\nu} (\gamma, t) - \frac{t^{-2\nu}}{\Gamma \l 1-2\nu \r} \, = \, \frac{^C\partial^{2\nu}}{\partial t^{2\nu}} \widetilde{\mathpzc{l}_\nu} (\gamma, t) \\
	& \frac{\partial^\nu}{\partial t^\nu} \widetilde{\mathpzc{l}_\nu} (\gamma, t) - \frac{t^{-\nu}}{\Gamma \l 1-\nu \r} \, = \, \frac{^C\partial^{\nu}}{\partial t^{\nu}} \widetilde{\mathpzc{l}_\nu} (\gamma, t) 
	\label{}
\end{align}
we therefore need to show that
\begin{equation} 
	\l \frac{^C\partial^{2\nu}}{\partial t^{2\nu}} + 2\lambda \frac{^C\partial^\nu}{\partial t^\nu} \r \widetilde{\mathpzc{l}_\nu} (\gamma, t) \, = \, - \gamma \widetilde{\mathpzc{l}_\nu} (\gamma, t).
	\label{}
\end{equation}
In light of
\begin{align}
	\frac{^{C}\partial^\nu}{\partial t^\nu} E_{\nu ,1} \l r_j t^\nu \r \, = \, r_j E_{\nu ,1} \l r_j t^\nu \r, \qquad j \, = \, 1, 2, 
	\label{}
\end{align}
\begin{align}
	 \frac{^{C}\partial^{2\nu}}{\partial t^{2\nu}} E_{\nu ,1} \l r_j t^\nu \r \, = \, r_j^2 \, E_{\nu, 1} \l r_j t^\nu \r + \frac{t^{-\nu} r_j}{\Gamma \l 1-\nu \r},
	\label{duestelle}
\end{align}
we are able to show that \eqref{laplace di lbold} solves \eqref{stellina}. We first check result \eqref{duestelle} as follows, for $0 < 2\nu < 1$
\begin{align}
	&  \frac{^{C}\partial^{2\nu}}{\partial t^{2\nu}} E_{\nu ,1} \l r_j t^{\nu } \r \,  =  \,  \sum_{k=0}^\infty \frac{r_j^k}{\Gamma \l \nu k + 1 \r} \frac{^{C}\partial^{2\nu}}{\partial t^{2\nu}} t^{\nu k} \notag \\
	 = \, & \sum_{k=1}^\infty \frac{r_j^k}{\Gamma \l \nu k + 1 \r} \frac{\nu k}{\Gamma \l 1-2\nu \r} \int_0^t s^{\nu k -1} \l t-s \r^{-2\nu} \, ds \notag \\
	 = \, & \sum_{k=1}^\infty \frac{r_j^k \, t^{\nu k - 2\nu}}{\Gamma \l \nu k  \r} \frac{1}{\Gamma \l 1-2\nu \r} \int_0^1 s^{\nu k -1} \l 1-s \r^{1-2\nu -1} \, ds \notag \\
	 = \, & \sum_{k=1}^\infty \frac{r_j^k \, t^{\nu k - 2 \nu}}{\Gamma \l \nu k - 2 \nu +1 \r} \, = \, \sum_{k=0}^\infty \frac{r_j^{k+1} \, t^{\nu k - \nu}}{\Gamma \l \nu k - \nu +1 \r} \notag \\
	 = \, & r_j t^{-\nu} \left[ \sum_{k=1}^\infty \frac{(r_j t^\nu)^k}{\Gamma \l \nu k - \nu + 1 \r} + \frac{1}{\Gamma \l 1-\nu \r} \right] \, = \, r_j^2 \, E_{\nu, 1} \l r_j t^\nu \r + \frac{t^{-\nu} r_j}{\Gamma \l 1-\nu \r}.
	\label{}
\end{align}

Therefore
\begin{align}
& \l \frac{^{C}\partial^{2\nu}}{\partial t^{2\nu}}  +2\lambda \frac{^{C}\partial^\nu}{\partial t^\nu} \r 	\widetilde{\mathpzc{l}_\nu} \l \gamma, t \r  \, = \notag \\
  = \, & \frac{1}{2} \left[ \l 1+ \frac{\lambda}{\sqrt{\lambda^2 - \gamma}} \r \frac{^{C}\partial^{2\nu}}{\partial t^{2\nu}} E_{\nu, 1} \l r_1 t^\nu \r + \l 1- \frac{\lambda}{\sqrt{\lambda^2 - \gamma}} \r \frac{^{C}\partial^{2\nu}}{\partial t^{2\nu}} E_{\nu,1} \l r_2 t^\nu \r \right] \,  \notag \\
&  + 2\lambda \frac{1}{2} \left[ \l 1+ \frac{\lambda}{\sqrt{\lambda^2 - \gamma}} \r \frac{^{C}\partial^{\nu}}{\partial t^{\nu}} E_{\nu, 1} \l r_1 t^\nu \r + \l 1- \frac{\lambda}{\sqrt{\lambda^2 - \gamma}} \r \frac{^{C}\partial^{\nu}}{\partial t^{\nu}} E_{\nu,1} \l r_2 t^\nu \r \right] \notag \\
= \, & \frac{1}{2} \left[ \l 1+ \frac{\lambda}{\sqrt{\lambda^2 - \gamma}} \r \l r_1^2 \, E_{\nu, 1} \l r_1 t^\nu \r + \frac{t^{-\nu} r_1}{\Gamma \l 1-\nu \r} \r  \right. \notag \\
 &  \left. + \l 1- \frac{\lambda}{\sqrt{\lambda^2 - \gamma}} \r \l r_2^2 \, E_{\nu, 1} \l r_2 t^\nu \r + \frac{t^{-\nu} r_2}{\Gamma \l 1-\nu \r} \r \right]  \notag \\
&  + 2\lambda \frac{1}{2} \left[ \l 1+ \frac{\lambda}{\sqrt{\lambda^2 - \gamma}} \r \l r_1 E_{\nu ,1} \l r_1 t^\nu \r \r + \l 1- \frac{\lambda}{\sqrt{\lambda^2 - \gamma}} \r r_2 E_{\nu ,1} \l r_2 t^\nu \r \right] \notag \\
 = \, & \frac{1}{2} \left[ r_1 \l 1+ \frac{\lambda}{\sqrt{\lambda^2 - \gamma}}  \r E_{\nu,1} \l r_1 t^\nu \r \l r_1 + 2\lambda \r + r_2 \l 1-\frac{\lambda}{\sqrt{\lambda^2 - \gamma}} \r \, \bm{\cdot} \right. \notag \\ 
 &  \bm{\cdot} E_{\nu,1} \l r_2 t^\nu \r \l r_2 + 2\lambda \r \bigg] \notag \\
 = \, & -\frac{\gamma}{2}   \frac{\lambda + \sqrt{\lambda^2 - \gamma}}{\sqrt{\lambda^2 - \gamma}}  E_{\nu, 1} \l r_1 t^\nu \r - \frac{\gamma}{2}   \frac{\sqrt{\lambda^2 - \gamma} -\lambda  }{\sqrt{\lambda^2 - \gamma}}  E_{\nu, 1} \l r_2 t^\nu \r \notag \\
 = \, & -\gamma \left[ \frac{1}{2} \left[ \l 1+ \frac{\lambda}{\sqrt{\lambda^2 - \gamma}} \r E_{\nu, 1} \l r_1 t^\nu \r + \l 1- \frac{\lambda}{\sqrt{\lambda^2 - \gamma}} \r E_{\nu, 1} \l r_2 t^\nu \r   \right] \right] \notag \\
 = \, & -\gamma \, \widetilde{\mathpzc{l}_\nu} \l \gamma, t \r.
\label{}
\end{align}
In the last steps we used the fact that
\begin{equation}
	\l 1+ \frac{\lambda }{\sqrt{\lambda^2 - \gamma}} \r \frac{r_1 \, t^{-\nu}}{\Gamma \l 1- \nu \r}+	\l 1- \frac{\lambda }{\sqrt{\lambda^2 - \gamma}} \r \frac{r_2 \, t^{-\nu}}{\Gamma \l 1- \nu \r} \, = \, 0,
	\label{}
\end{equation}
and
\begin{equation}
r_1 + 2\lambda \, = \, -r_2, \qquad r_2 + 2\lambda \, = \, -r_1, \qquad r_1r_2 \, = \, \gamma.
\label{}
\end{equation}
\end{proof}
\end{te}

\begin{os} \normalfont
The derivation of result \eqref{laplace di lbold} suggests an alternative proof for the Fourier transform (Theorem 2.2 in \citet{ptrf}) of the law of the time-fractional telegraph process.
\end{os}

\begin{os} \normalfont
From \eqref{double laplace di lnu} we get the time Laplace transform of $\mathpzc{l}_\nu (x, t)$, for $x>0, \mu > 0, 0 < \nu < \frac{1}{2} $, as
\begin{equation} 
	\widetilde{\mathpzc{l}_\nu} \l x, \mu  \r \, = \, \mu^{2\nu -1} e^{-x\mu^{2\nu}} e^{-2\lambda x \mu^\nu}  + 2\lambda \mu^{\nu -1} e^{-2\lambda x \mu^\nu} e^{-x\mu^{2\nu}}.
	\label{questaqui}
\end{equation}
Since (see formulas \eqref{155} and \eqref{lapl l})
\begin{equation} 
	\widetilde{l_\nu} (x, \mu) \, = \, \int_0^\infty e^{-\mu t} \frac{1}{t^\nu} W_{-\nu, 1-\nu} \l -\frac{x}{t^\nu} \r \, dt	\, = \, \mu^{\nu -1} e^{-x\mu^\nu}
	\label{}
\end{equation}
and (see formula \eqref{153})
\begin{equation} 
	\widetilde{h_\nu} (\mu, t) \, = \, \int_0^\infty e^{-\mu x} h_\nu (x, t) \, dx \, = \, e^{-t\mu^\nu},
	\label{}
\end{equation}
we are able to invert \eqref{questaqui} and we obtain the explicit distribution of the process $\mathpzc{L}^\nu (t)$, $t>0$, which reads
\begin{align}
	\mathpzc{l}_\nu (x, t) \,  = \, & \frac{\Pr \ll \mathpzc{L}^\nu (t) \in dx \rr}{dx} \, \notag \\
	 = \, & \int_0^t l_{2\nu} \l  x, s \r \, h_\nu (t-s, 2\lambda x) \, ds + 2\lambda \int_0^t l_\nu (2\lambda x, s) \, h_{2\nu} \l t-s, x \r \, ds  \notag \\ 
	 = \, & \int_0^t \frac{1}{s^{2\nu}} \, W_{-2\nu, 1-2\nu} \l -\frac{x}{s^{2\nu}} \r \,  h_\nu (t-s, 2\lambda x )  \, ds   \notag \\
	 & + 2\lambda \int_0^t \frac{1}{s^\nu} \, W_{-\nu, 1-\nu} \l -\frac{2\lambda x}{s^\nu} \r \, h_{2\nu} \l t-s, x \r \, ds.
	\label{legge esplicita di lbold per ogni nu}
\end{align}
The densities $h_\nu$ and $h_{2\nu}$ can be written down in terms of series expansion of stable laws (see pag. 245 of \citet{orsann}).
\end{os}

\section{$n$-dimensional stable laws and fractional Laplacian}
\label{laplaciansec}
Let 
\begin{equation} 
\bm{S}_n^{2\beta}(t) = \l S_1^{2\beta}(t), S_2^{2\beta}(t), \cdots , S_n^{2\beta}(t) \r, \qquad t>0, \beta \in (0,1],
	\label{isotropic vector}
\end{equation}
be the isotropic stable $n$-dimensional process with joint characteristic function
\begin{align}
	\widehat{v_n^{2\beta}} \l \bm{\xi}, t \r \, & = \, \widehat{v_n^{2\beta}} \l \xi_1, \xi_2, \cdots, \xi_n, t \r \, = \, \mathbb{E}e^{i \bm{\xi} \cdot \bm{S}_n^{2\beta}(t)} \, = \, e^{-t \l \sqrt{\xi_1^2 + \xi_2^2 + \cdots + \xi_n^2} \r^{2\beta}} \notag \\
	& = \, e^{-t \left\| \bm{\xi} \right\|^{2\beta}}.
	\label{caratteristica stable ndim}
\end{align}
The density corresponding to the characteristic function $\widehat{v_n^{2\beta}} \l \bm{\xi}, t \r$ is given by
\begin{equation} 
	v_n^{2\beta} \l \bm{x}, t \r \, = \, v_n^{2\beta} \l x_1, x_2, \cdots , x_n , t  \r \, = \, \frac{1}{\l 2\pi \r^n} \int_{\mathbb{R}^n}  e^{-i \bm{\xi} \cdot \bm{x}} e^{-t \left\| \bm{\xi} \right\|^{2\beta}} d\bm{\xi}.
	\label{}
\end{equation}
The equation governing the distribution $v_n^{2\beta} \l \bm{x}, t \r$ of the vector process $\bm{S}_n^{2\beta}(t)$, $t>0$, is
\begin{equation} 
	 \l \frac{\partial}{\partial t}  + \l - \Delta \r^{\beta} \r v_n^{2\beta} \l \bm{x}, t \r \, = \, 0, \qquad \bm{x} \in \mathbb{R}^n, t>0,
	\label{}
\end{equation}
where the fractional negative Laplacian is related to the classical Laplacian by means of the following relationships (Bochner representation, see for example \citet{balakrishnan, bochner})
\begin{align}
	 &\frac{\sin \pi \beta}{\pi} \int_0^\infty d\lambda \, \lambda^{\beta -1} \l \lambda - \Delta \r^{-1} \, \Delta \, = \, \frac{\sin \pi \beta}{\pi} \int_0^\infty \lambda^{\beta -1} \l \int_0^\infty e^{-w \l \lambda - \Delta \r} dw \r \Delta \, d\lambda \notag \\
	 = \, &\frac{\sin \pi \beta}{\pi} \, \Delta \, \Gamma (\beta) \int_0^\infty w^{1-\beta -1} e^{-w \l -\Delta \r} dw \, = \, \frac{\Delta}{\Gamma \l 1-\beta \r}\int_0^\infty w^{1-\beta -1} e^{-w \l - \Delta \r} dw \notag \\
	 = \, & - \l - \Delta \r^\beta.
	\label{}
\end{align}
A definition of the fractional negative Laplacian can be given in the space of the Fourier transforms as follows
\begin{align}
	- \l -\Delta \r^\beta u(\bm{x}) \, = \, - \frac{1}{(2\pi)^n} \int_{\mathbb{R}^n} e^{-i \bm{x} \cdot \bm{\xi}} \l \xi_1^2 + \xi_2^2 + \cdots + \xi_n^2 \r^\beta \, \widehat{u} \l \bm{\xi} \r \, d\bm{\xi},
	\label{def laplacian ft}
\end{align}
where
\begin{align}
	\textrm{Dom} \l - \Delta \r^\beta \, = \, \left\lbrace u \in L_{ \textrm{loc} }^1 \l \mathbb{R}^n \r \, : \, \int_{\mathbb{R}^n} \left| \widehat{u} \l \bm{\xi} \r \right|^2 \l 1+ \left\| \bm{\xi}  \right\|^{2\beta} \r \, d\bm{\xi} < \infty \right\rbrace.
	\label{}
\end{align}

An equivalent alternative definition of the $n$-dimensional fractional Laplacian is
\begin{equation} 
	\l - \Delta \r^\beta u(\bm{x}) \, = \, c(\beta, n) \, \textrm{P.V.}  \int_{\mathbb{R}^n} \frac{u(\bm{x})-u(\bm{y})}{\left\| \bm{x}-\bm{y} \right\|^{n+2\beta}} dy,
	\label{singolar}
\end{equation} 
where the multiplicative constant $c(\beta, n)$ must be evaluated in such a way that
\begin{align}
	\int_{\mathbb{R}^n} e^{i\bm{\xi} \cdot \bm{x}} \l -\Delta \r^\beta \, u(\bm{x}) \, d\bm{x} \, = \,  \left\| \bm{\xi} \right\|^{2\beta}  \int_{\mathbb{R}^n} e^{i \bm{\xi} \cdot \bm{x}} u(\bm{x}) \, d\bm{x}.
	\label{}
\end{align}
Let us focus our attention on the one-dimensional case of \eqref{singolar}. In this case we have that, for $0<2\beta<1$,
\begin{align}
	&  \l - \frac{\partial^2}{\partial x^2} \r^\beta  u(x) \,  = \,  c(\beta, 1) \, \textrm{P.V.} \int_{\mathbb{R}} \frac{u(x) - u(y)}{\left| x-y \right|^{1+2\beta}} \notag \\	
	 = \, & c(\beta, 1) \, \lim_{\epsilon \to 0} \left[ \int_{-\infty}^{0-\epsilon} \frac{u(x)-u(x-z)}{|z|^{1+2\beta}} dz + \int_{0+\epsilon}^\infty \frac{u(x) - u(x-z)}{|z|^{1+2\beta}} dz \right] \notag \\
	 = \, & c(\beta, 1) \, \lim_{\epsilon \to 0} \left[ \int_{0+\epsilon}^\infty \frac{u(x)-u(x+z)}{z^{1+2\beta}} dz + \int_{0+\epsilon}^\infty \frac{u(x) - u(x-z)}{z^{1+2\beta}} dz \right] \notag \\
	 = \, & \frac{\Gamma (1-2\beta)}{2\beta} c(\beta, 1) \left[ \frac{1}{\Gamma \l 1-2\beta \r} \frac{d}{dx} \l \int_{-\infty}^x \frac{u(z)\; dz}{(x-z)^{2\beta}}  -   \int_x^\infty \frac{u(z)\; dz}{(z-x)^{2\beta}} \r \right],
	\label{}
\end{align}
where in the intermediate steps, we considered the relation between the Marchaud and the Weyl fractional derivatives. By setting
\begin{equation} 
	c(\beta, 1) \, = \, \frac{2\beta}{2\, \Gamma \l 1-2\beta \r \; \cos \beta \pi},
	\label{}
\end{equation}
we have that, for $0<2\beta<1$,
\begin{align}
	&  -\l -\frac{\partial^2}{\partial x^2} \r^\beta  u(x) \,  = \notag \\
	 = \, & -\frac{1}{2\cos \beta \pi} \left[\frac{1}{\Gamma (1-2\beta)} \frac{d}{dx} \int_{-\infty}^x \frac{u(z)\; dz}{(x-z)^{2\beta}} - \frac{1}{\Gamma (1-2\beta)} \frac{d}{dx} \int_x^\infty \frac{u(z)\; dz}{(z-x)^{2\beta}} \right] \notag \\
	 = \, &  -\frac{1}{2\cos \beta \pi} \frac{1}{\Gamma (1-2\beta)} \frac{d}{dx} \int_{-\infty}^\infty \frac{u(z)}{|x-z|^{2\beta}} dz \, = \, \frac{\partial^{2\beta}}{\partial |x|^{2\beta}} u(x),
	\label{312}
\end{align}
where $	\frac{\partial^{2\beta}}{\partial |x|^{2\beta}}$ represents the Riesz operator.

\begin{os} \normalfont
We notice that, for $0<2\beta < 1$,
\begin{equation} 
	\mathcal{F} \left[ \frac{\partial^{2\beta}}{\partial |x|^{2\beta}} u(x) \right] (\xi) \, = \, -| \xi |^{2\beta} \, \widehat{u}(\xi).
	\label{3113}
\end{equation}
This is due to the calculation
\begin{align}
		&  \mathcal{F} \left[ \frac{\partial^{2\beta}}{\partial |x|^{2\beta}} u(x) \right] (\xi) \, = \notag \\
		 = \, &  -\frac{1}{2\cos \beta \pi} \frac{1}{\Gamma \l 1-2\beta \r} \left[  \int_{-\infty}^\infty dx \, e^{i\xi x} \l \frac{d}{dx} \int_{-\infty}^x \frac{u(z) \; dz}{\l x-z \r^{2\beta}} - \frac{d}{dx} \int_x^\infty \frac{u(z) \; dz}{\l z-x \r^{2\beta}} \,  \r \right] \notag \\
		 = \, & \frac{i \xi}{2\cos \beta \pi} \frac{1}{\Gamma \l 1-2\beta \r} \left[  \int_{-\infty}^\infty dx \, e^{i\xi x} \l  \int_{-\infty}^x \frac{u(z) \; dz}{\l x-z \r^{2\beta}} -  \int_x^\infty \frac{u(z) \; dz}{\l z-x \r^{2\beta}} \,  \r \right] \notag \\
		 = \, & \frac{i \xi}{2\cos \beta \pi} \frac{1}{\Gamma \l 1-2\beta \r} \left[ \int_{-\infty}^\infty dz \, u(z) \l  \int_{z}^\infty \frac{e^{i\xi x} \, dx}{\l x-z \r^{2\beta}} -  \int_{-\infty}^z \frac{e^{i\xi x} \, dx}{\l z-x \r^{2\beta}} \,  \r \right] \notag \\
		 = \, & \frac{i\xi}{2\cos \beta \pi} \frac{1}{\Gamma \l 1-2\beta \r} \left[ \int_{-\infty}^\infty  e^{i\xi z} u(z) \, dz \l \int_0^\infty \frac{e^{i\xi y} }{y^{2\beta}} \, dy - \int_0^\infty \frac{e^{-i\xi y}}{y^{2\beta}} \, dy \r \right] \notag \\
		 = \, & -\frac {2\xi}{2\cos \beta \pi} \frac{1}{\Gamma \l 1-2\beta \r} \int_{-\infty}^\infty e^{i\xi z} u(z) \, dz \int_0^\infty \frac{\sin \xi y}{y^{2\beta}} \, dy \notag \\
		 = \, & -\frac {\xi}{\cos \beta \pi} \frac{1}{\Gamma \l 1-2\beta \r} \;  \frac{\widehat{u}(\xi)}{\Gamma \l 2\beta \r}  \int_0^\infty \int_0^\infty \sin \xi y \,  e^{-wy} w^{2\beta -1} \; dw \, dy \notag \\
		 = \, & -\frac {\xi}{\cos \beta \pi} \frac{1}{\Gamma \l 1-2\beta \r} \;  \frac{\widehat{u}(\xi)}{\Gamma \l 2\beta \r}  \int_0^\infty dw \, w^{2\beta -1} \int_0^\infty dy \, e^{-wy} \l \frac{e^{i\xi y}-e^{-i\xi y}}{2i} \r \notag \\
		 = \, & -\frac {\xi^2}{\cos \beta \pi} \frac{1}{\Gamma \l 1-2\beta \r} \;  \frac{\widehat{u}(\xi)}{\Gamma \l 2\beta \r}   \int_0^\infty dw \frac{w^{2\beta -1}}{w^2 + \xi^2} \notag \\
		 = \, & -\frac {\xi^2}{\cos \beta \pi} \frac{1}{\Gamma \l 1-2\beta \r} \;  \frac{\widehat{u}(\xi)}{\Gamma \l 2\beta \r}   \int_0^\infty dw \, w^{2\beta -1} \int_0^\infty dy \, e^{-y \l w^2 + \xi^2 \r} \notag \\
		 = \, &  -\frac {\xi^2}{2\cos \beta \pi} \frac{1}{\Gamma \l 1-2\beta \r} \;  \frac{\widehat{u}(\xi)}{\Gamma \l 2\beta \r} \frac{\Gamma \l \beta \r \, \Gamma \l 1-\beta \r}{| \xi |^{2-2\beta}} \, = \, -  | \xi |^{2\beta} \; \widehat{u}(\xi).
	\label{}
\end{align}
This concludes the proof of \eqref{3113}.
\end{os}

\section{Space-time fractional telegraph equation}

We consider now the composition of an isotropic vector of stable processes $\bm{S}_n^{2\beta} (t)$, $t>0$, defined in \eqref{isotropic vector}, with the positively-valued process, defined in \eqref{construction L},
\begin{equation} 
	\mathpzc{L}^\nu (t) \, = \, \inf \ll s>0: \mathpzc{H}^\nu (s) \, = \, H_1^{2\nu} (s) + (2\lambda )^{ \frac{1}{\nu} } H_2^\nu (s) \geq t  \rr, \qquad t>0,
	\label{}
\end{equation}
where $H_1^{2\nu}$, $H_2^{\nu}$ are independent positively skewed stable processes of order $2\nu$ and $\nu$, respectively.
The distribution $w_\nu^\beta \l \bm{x}, t \r$ of the process $\bm{S}_n^{2\beta} \l c^2 \mathpzc{L}^\nu (t) \r$, $t>0$, $\beta \in (0,1]$, is the fundamental solution to the space-time fractional telegraph equation
\begin{equation} 
	\l \frac{^C\partial^{2\nu}}{\partial t^{2\nu}} + 2\lambda \frac{^C\partial^\nu}{\partial t^\nu} \r w_\nu^\beta \l \bm{x}, t \r \, = \, -c^2 \l -\Delta \r^\beta \, w_\nu^\beta \l \bm{x}, t \r, \qquad \bm{x} \in \mathbb{R}^n, t>0.
	\label{stfte}
\end{equation}
In our view the next theorem generalizes some previous results because we here have fractionality in space and time and the equation \eqref{stfte} is defined in $\mathbb{R}^n$.

\begin{te}
\label{theorem composition}
For $\nu \in \l 0, \frac{1}{2} \right]$, $\beta \in (0,1]$ and $c>0$ the solution to the Cauchy problem for the space-time fractional $n$-dimensional telegraph equation
\begin{align}
	\begin{cases}
	\l \frac{^C\partial^{2\nu}}{\partial t^{2\nu}} + 2\lambda \frac{^C\partial^\nu}{\partial t^\nu} \r w_\nu^\beta \l \bm{x}, t \r \, = \, -c^2 \l -\Delta \r^\beta \, w_\nu^\beta \l \bm{x}, t \r, \qquad \bm{x} \in \mathbb{R}^n, t>0 \\
	w_\nu^\beta \l \bm{x}, 0 \r \, = \, \delta \l \bm{x} \r,
	\end{cases}
	\label{stftcp}
\end{align}
coincides with the probability law of the vector process
\begin{equation} 
	\bm{W}_n (t) \, = \, \bm{S}_n^{2\beta} \l c^2 \mathpzc{L}^\nu (t) \r, \qquad t>0,
	\label{composizione}
\end{equation}
and has Fourier transform which reads
\begin{align}
	&\widehat{w_\nu^\beta} \l \bm{\xi}, t \r \, = \notag \\
	= \, & \frac{1}{2} \left[ \l 1+ \frac{\lambda}{\sqrt{\lambda^2- c^2 \left\| \bm{\xi} \right\|^{2\beta}}} \r E_{\nu,1} \l r_1 t^\nu \r + \l 1-\frac{\lambda}{\sqrt{\lambda^2 - c^2 \left\| \bm{\xi} \right\|^{2\beta}}} \r E_{\nu,1} \l r_2 t^\nu \r  \right],
	\label{fourier mittagleffler}
\end{align}
where 
\begin{equation} 
	r_1 \, = \, -\lambda + \sqrt{\lambda^2 - c^2 \left\| \bm{\xi} \right\|^{2\beta}}, \qquad r_2 \, = \, -\lambda -\sqrt{\lambda^2 - c^2 \left\| \bm{\xi} \right\|^{2\beta}}.
	\label{}
\end{equation}
The time derivatives appearing in \eqref{stftcp} must be meant in the Dzerbayshan-Caputo sense. The fractional Laplacian is defined in \eqref{def laplacian ft}.
\begin{proof}
By taking the Laplace transform of \eqref{stftcp} we have
\begin{align}
	\mu^{2\nu} \widetilde{w_\nu^\beta } \l \bm{x}, \mu \r - \mu^{2\nu -1} \delta (\bm{x}) + 2\lambda \left[ \mu^\nu \widetilde{w_\nu^\beta} \l \bm{x}, \mu \r - \mu^{\nu -1} \delta (\bm{x}) \right] \, = \, - c^2 \l -\Delta \r^\beta \widetilde{w_\nu^\beta} \l \bm{x}, \mu \r,
	\label{laplace stftcp}
\end{align}
where we used the fact that (see \cite{kill} page 98, Lemma 2.24)
\begin{equation} 
	\mathcal{L} \left[ \frac{^C\partial^\nu}{\partial t^\nu} w_\nu^\beta \l \bm{x}, t \r \right] \, = \, \mu^\nu \widetilde{w_\nu^\beta} (\bm{x}, \mu) - \mu^{\nu -1} w_\nu^\beta (\bm{x}, 0).
	\label{}
\end{equation}
Now the Fourier transform of \eqref{laplace stftcp} yields
\begin{equation} 
	\l \mu^{2\nu} + 2\lambda \mu^\nu \r \widehat{\widetilde{w_\nu^\beta}} \l \bm{\xi}, \mu \r - \l \mu^{2\nu -1} + 2\lambda \mu^{\nu -1 } \r \, = \, -c^2  \left\| \bm{\xi}  \right\|^{2\beta} \widehat{\widetilde{w_\nu^\beta}} \l \bm{\xi}, \mu \r,
	\label{}
\end{equation}
and thus
\begin{equation} 
	\widehat{\widetilde{w_\nu^\beta}} \l \bm{\xi}, \mu \r \, = \, \frac{\mu^{2\nu -1} + 2\lambda \mu^{\nu -1}}{\mu^{2\nu} + 2\lambda \mu^\nu + c^2 \left\| \bm{\xi} \right\|^{2\beta}}, \qquad \mu >0, \bm{\xi} \in \mathbb{R}^n.
	\label{hat tilde}
\end{equation}
The probability density of the process $\bm{W}_n(t)$, $t>0$, defined in \eqref{composizione}, can be written as
\begin{equation} 
	w_\nu^\beta \l \bm{x}, t \r  \, = \, \int_0^\infty v_\beta \l \bm{x}, c^2 s \r \, \mathpzc{l}_\nu \l s, t \r \, ds,
	\label{}
\end{equation}
and has Fourier transform equal to
\begin{equation} 
	\int_{\mathbb{R}^n} e^{i \bm{\xi} \cdot \bm{x}} \, w_\nu^\beta \l \bm{x}, t \r \, d\bm{x} \, = \, \int_0^\infty e^{-c^2s\left\| \bm{\xi} \right\|^{2\beta}} \, \mathpzc{l}_\nu (s, t) \, ds.
	\label{fourier tbold}
\end{equation}
In order to show that the Laplace transform of \eqref{fourier tbold} concides with \eqref{hat tilde}, we have to derive the Laplace transform of $\mathpzc{l}_\nu (x, t)$, with respect to the time $t$. Since
\begin{equation} 
	\Pr \ll \mathpzc{L}^\nu (t) < x  \rr \, = \, \Pr \ll \mathpzc{H}^\nu (x) >t \rr
	\label{}
\end{equation}
we have that
\begin{align}
	&\widetilde{\mathpzc{l}_\nu} \l x, \mu \r  =  \notag \\
	 = \, &\int_0^\infty e^{-\mu t} \frac{\partial}{\partial x} \int_t^\infty \Pr \ll \mathpzc{H}^\nu (x) \in ds \rr \, dt \,  = \,  \int_0^\infty e^{-\mu t} \l - \frac{\partial}{\partial x} \int_0^t \mathpzc{h}_\nu \l s, x \r  ds \r  dt \notag \\
	 = \, &-\frac{\partial}{\partial x} \frac{e^{-x \mu^{2\nu}-2\lambda x \mu^\nu}}{\mu} \, = \, \l \mu^{2\nu -1} + 2\lambda \mu^{\nu -1} \r \, e^{-x\mu^{2\nu}-2\lambda x \mu^\nu},
	\label{laplace lbold inv}
\end{align}
where we used result \eqref{lapl proc}.
Now we can complete the proof by taking the Laplace transform of \eqref{fourier tbold} so that, in view of \eqref{laplace lbold inv}, we obtain
\begin{align}
	 & \int_0^\infty e^{-\mu t} dt \int_0^\infty e^{-c^2s\left\| \bm{\xi} \right\|^{2\beta}} \, \mathpzc{l}_\nu (s, t) \, ds \, = \notag \\
	  = \, &  \l \mu^{2\nu -1} + 2\lambda \mu^{\nu -1} \r \, \int_0^\infty e^{-sc^2\left\| \bm{\xi} \right\|^{2\beta}-s\mu^{2\nu}-2\lambda s \mu^\nu} \, ds \, = \, \frac{\mu^{2\nu -1} + 2\lambda \mu^{\nu -1}}{\mu^{2\nu} + 2\lambda \mu^\nu + c^2 \left\| \bm{\xi} \right\|^{2\beta}},
	\label{fourier laplace caso generale}
\end{align}
which coincides with \eqref{hat tilde}. The unicity of Fourier-Laplace transform proves that the claimed result holds. The proof that the Fourier transform of $w_\nu^\beta \l \bm{x}, t \r$ has the form \eqref{fourier mittagleffler} can be carried out by means of the calculation performed in Theorem \ref{teorema inverso}. We have that
\begin{align}
	\widehat{\widetilde{w_\nu^\beta}} \l \bm{\xi}, \mu \r \, & = \, \frac{\mu^{2\nu-1 }+2\lambda \mu^{\nu-1}}{\mu^{2\nu}+ 2\lambda \mu^\nu  + c^2 \left\| \bm{\xi} \right\|^{2\beta} } \, = \, \frac{\mu^{\nu-1}}{\mu^\nu - r_1} + \frac{\mu^{\nu-1}}{\mu^\nu - r_2} - \frac{\mu^{2\nu-1}}{\l \mu^\nu - r_1 \r \l \mu^\nu - r_2 \r} \notag \\
	& = \, \frac{\mu^{\nu-1}}{\mu^\nu - r_1} + \frac{\mu^{\nu-1}}{\mu^\nu - r_2} - \left[ \frac{\mu^{\nu- (1-\nu)}}{\mu^\nu - r_1} - \frac{\mu^{\nu-(1-\nu)}}{\mu^\nu - r_2} \right] \, \frac{1}{2\sqrt{\lambda^2 - c^2 \left\| \bm{\xi} \right\|^{2\beta}}},
	\label{byinverting}
\end{align}
where
\begin{equation} 
	r_1 \, = \, -\lambda + \sqrt{\lambda^2 - c^2 \left\| \bm{\xi} \right\|^{2\beta}}, \qquad r_2 \, = \, -\lambda -\sqrt{\lambda^2 - c^2 \left\| \bm{\xi} \right\|^{2\beta}}.
	\label{}
\end{equation}
and thus by inverting \eqref{byinverting} by means of \eqref{proprietàinversione}, we obtain result \eqref{fourier mittagleffler}. An alternative derivation of \eqref{fourier mittagleffler} can be carried out as follows
\begin{align}
	\widehat{w_\nu^\beta} \l \bm{\xi}, t \r \, & = \, \int_{-\infty}^\infty e^{i \bm{\xi} \cdot \bm{x}} d\bm{x} \int_0^\infty \Pr \ll \bm{S}_n^{2\beta} \l c^2 s \r \in d\bm{x} \rr \Pr \ll \mathpzc{L}^\nu (t) \in ds \rr \notag \\
	& = \, \int_0^\infty e^{-c^2 s \left\| \bm{\xi} \right\|^{2\beta} } \Pr \ll \mathpzc{L}^\nu (t) \in ds \rr \, = \, \eqref{fourier mittagleffler}
	\label{}
\end{align}
because of Theorem \ref{teorema inverso}.
\end{proof}
\end{te}

\subsection{The case $\nu = \frac{1}{2}$, subordinator with drift}

The  fractional equation \eqref{stfte}, for $n=1$, $\nu = \frac{1}{2}$, reads
\begin{equation} 
	\l \frac{\partial}{\partial t} + 2\lambda \frac{^C\partial^{ \frac{1}{2} }}{\partial t^{ \frac{1}{2} }} \r \, w_{ \frac{1}{2} }^\beta (x, t) \, = \, c^2 \l \frac{\partial^{2\beta}}{\partial |x|^{2\beta}} \r \, w_{ \frac{1}{2} }^\beta \l x, t \r, \qquad 0 < \beta < 1,
	\label{51}
\end{equation}
where $\frac{\partial^{2\beta}}{\partial |x|^{2\beta}}$ is the Riesz operator defined in \eqref{312}. For $\beta = 1$ we have the special case
\begin{equation} 
	\l \frac{\partial}{\partial t} + 2\lambda \frac{^C\partial^{ \frac{1}{2} }}{\partial t^{ \frac{1}{2} }} \r \, w_{ \frac{1}{2} }^1 \l x, t \r \, = \, c^2 \frac{\partial^2}{\partial x^2} \, w_{ \frac{1}{2} }^1 \l x, t \r
	\label{}
\end{equation}
dealt with in \citet{ptrf}. The construction of the composition related to equation \eqref{51} involves the subordinator
\begin{equation} 
	\mathpzc{H}^{ \frac{1}{2} } (t) \, = \, t + \l 2\lambda \r^2 H^{ \frac{1}{2} } (t), \qquad t>0,
	\label{H un mezzo}
\end{equation}
where $H^{ \frac{1}{2} } (t)$, $t>0$, is a positively-skewed stable process and has the same law as the first-passage time of a Brownian motion through level $ \frac{t}{\sqrt{2}} $. We note that $\mathpzc{H}^{ \frac{1}{2} } (t)$, $t>0$, has distribution with support $[t, \infty )$ and thus differs from $\mathpzc{H}^\nu (t)$, $t>0$, $0< \nu < \frac{1}{2}$, which instead has support $[0, \infty)$.
The distribution of \eqref{H un mezzo} writes
\begin{equation} 
	\Pr \ll \mathpzc{H}^{ \frac{1}{2} } (t)  < x \rr \, = \, \int_{ 0 }^{\frac{x-t}{\l 2\lambda \r^2}} \frac{t}{\sqrt{2}} \frac{e^{-\frac{t^2}{4z}}}{\sqrt{2\pi z^3}} \, dz, \qquad x>t>0.
	\label{}
\end{equation}
The inverse process
\begin{equation} 
	\mathpzc{L}^{ \frac{1}{2} } (t) \, = \, \inf \ll s:s+\l 2\lambda \r^2 H^{ \frac{1}{2} } (s) \geq t \rr \, = \, \inf \ll s:\mathpzc{H}^{ \frac{1}{2} } (s) \geq t \rr
	\label{L un mezzo}
\end{equation}
is related to \eqref{H un mezzo} by means of the relationship
\begin{equation} 
	\Pr \ll \mathpzc{L}^{ \frac{1}{2} } (t) < x \rr \, = \, \Pr \ll \mathpzc{H}^{ \frac{1}{2} } (x) > t \rr \, = \, \int_{ \frac{t-x}{\l 2\lambda \r^2} }^\infty \frac{x}{\sqrt{2}} \frac{e^{ -\frac{x^2}{4z} }}{\sqrt{2\pi z^3}} \, dz.
	\label{56}
\end{equation}
From \eqref{56} we can extract the distributon of $\mathpzc{L}^{ \frac{1}{2} } (t)$, $t>0$, in the following manner
\begin{align}
	\mathpzc{l}_{ \frac{1}{2} } (x, t) \, & = \, \frac{\Pr \ll \mathpzc{L}^{ \frac{1}{2} } (t) \in dx \rr}{dx} \, = \, \frac{\partial}{\partial x} \int_{ \frac{t-x}{\l 2\lambda \r^2} }^\infty \frac{x \, e^{-\frac{x^2}{4z}}}{\sqrt{4\pi z^3}} \, dz \notag \\
	& = \, \frac{2\lambda x \, e^{ -\frac{\l 2\lambda x \r^2}{4(t-x)} }}{\sqrt{4\pi \l t-x \r^3}} + 2\lambda \frac{e^{ -\frac{\l 2\lambda x \r^2}{4(t-x)} }}{\sqrt{\pi (t-x)}}, \qquad 0<x<t.
	\label{l un mezzo}
\end{align}
\begin{os} \normalfont
The distribution \eqref{l un mezzo} can be also obtained from the general case \eqref{legge esplicita di lbold per ogni nu} which for $\nu = \frac{1}{2}$ becomes, for $0<x<t$,
\begin{align}
	\mathpzc{l}_{ \frac{1}{2} } (x, t) \, & = \, \int_0^t \delta (s-x) \, h_{\frac{1}{2}} (t-s, 2\lambda x) \, ds +2\lambda \int_0^t  l_{\frac{1}{2}} (2\lambda x, s) \, \delta \l x-(t-s) \r \, ds \notag \\
	& = \, h_{ \frac{1}{2} } \l t-x, 2\lambda x \r + 2\lambda \, l_{ \frac{1}{2} } (2\lambda x, t-x) \notag \\
	& = \, \frac{2\lambda x \, e^{ -\frac{\l 2\lambda x \r^2}{4(t-x)} }}{\sqrt{4\pi \l t-x \r^3}} + 2\lambda \frac{e^{ -\frac{\l 2\lambda x \r^2}{4(t-x)} }}{\sqrt{\pi (t-x)}}.
	\label{}
\end{align}
In the last step we used the fact that
\begin{equation} 
	L^{\frac{1}{2}} (t) \, \stackrel{law}{=} \, \left| B(t) \right|, \qquad t>0,
	\label{}
\end{equation}
where $L^{\frac{1}{2}} (t)$, $t>0$, dealt with in section \ref{prelimin}, is the inverse of the totally positively-skewed stable process $H^{\frac{1}{2}} (t)$, $t>0$.
\end{os}

The $t$-Laplace transform of \eqref{l un mezzo} becomes
\begin{align}
	\widetilde{\mathpzc{l}_{ \frac{1}{2} }} (x, \mu) \, & = \,  \int_x^\infty e^{-\mu t} \, \mathpzc{l}_{ \frac{1}{2} } (x, t) \, dt \, = \notag \\
	& = \, \frac{2\lambda x	}{\sqrt{2}} \int_x^\infty e^{-\mu t} \frac{e^{ -\frac{\l 2\lambda x \r^2}{4(t-x)} }}{\sqrt{2\pi (t-x)^3}} \, dt + 2\lambda \int_x^\infty e^{-\mu t} \frac{e^{-\frac{\l 2\lambda x \r^2}{4(t-x)}}}{\sqrt{\pi (t-x)}} \, dt \notag \\
	& = \, \frac{2\lambda x}{\sqrt{2}} e^{-\mu x} \int_0^\infty e^{-\mu t} \frac{e^{ -\frac{\l 2\lambda x \r^2}{4t} }}{\sqrt{2\pi t^3}} dt + 2\lambda e^{-\mu x} \int_0^\infty e^{-\mu t} \frac{e^{-\frac{\l 2\lambda x \r^2}{4t}}}{\sqrt{\pi t}} \, dt \notag \\
	& = \, e^{-\mu x} e^{-2\lambda x \sqrt{\mu}} + 2\lambda \, \mu^{ -\frac{1}{2} } \, e^{-\mu x} e^{-2\lambda x \sqrt{\mu}}.
	\label{tlaplace l un mezzo}
\end{align}
Finally the $x$-Laplace transform of \eqref{tlaplace l un mezzo} becomes
\begin{align}
	\widetilde{\widetilde{\mathpzc{l}_{ \frac{1}{2} }}} \l \gamma, \mu \r \, & = \, \int_0^\infty e^{-\gamma x} \l \int_x^\infty e^{-\mu t} \, \mathpzc{l}_{ \frac{1}{2} } (x, t) \, dt \r \, dx \notag \\
	& = \, \frac{1}{\mu + \gamma + 2\lambda \sqrt{\mu}} + \frac{2\lambda}{\sqrt{\mu}} \frac{1}{\mu + \gamma + 2\lambda \sqrt{\mu}} \, = \, \frac{1+ 2\lambda \mu^{ -\frac{1}{2} }  }{\mu + \gamma + 2\lambda \sqrt{\mu}},
	\label{laplace di l un mezzo}
\end{align}
which coincides with \eqref{double laplace di lnu}, for $\nu = \frac{1}{2}$.
Let us now consider the process $\bm{W}_n (t) = \bm{S}_{n}^{2\beta} \l c^2 \mathpzc{L}^\nu (t) \r$, $t>0$, dealt with in Theorem \ref{theorem composition}. For $\beta =1$, $n=1$ and $\nu = \frac{1}{2}$ this process becomes
\begin{equation} 
	W_1 (t) \, = \, S_1^2 \l c^2 \mathpzc{L}^{ \frac{1}{2} } (t) \r \, = \, B \l c^2 \mathpzc{L}^{ \frac{1}{2} } (t) \r, \qquad t >0
	\label{W uno un mezzo}
\end{equation}
where $B$ represents a standard Brownian motion and $\mathpzc{L}^{ \frac{1}{2} } (t) $, $t>0$, is the process defined in \eqref{L un mezzo}. With
\begin{equation} 
	p_{|B|} (x, t) \, = \, \frac{e^{-\frac{x^2}{4t}}}{\sqrt{\pi t}}, \qquad x >0, t>0,
	\label{}
\end{equation}
we denote the law of the process $\left| B(t) \right|$, $t>0$. In view of the previous results we are able to prove the following theorem.
\begin{te}
\label{teorema un mezzo}
The law of \eqref{W uno un mezzo} coincides with the law of the composition
\begin{equation} 
	 \mathpzc{W} (t) \, = \, T \l \left| B(t) \right| \r, \qquad t >0,
	\label{}
\end{equation}
where $T$ is the telegraph process \eqref{costruzionetelegrafo} with parameters $c>0$, $\lambda >0$ and law $p_T (x, t)$ which has characteristic function
\begin{align}
	& \widehat{p_T} \l \xi, t \r \, = \notag \\
	 = \, & \frac{1}{2} \left[ \l 1+\frac{\lambda}{\sqrt{\lambda^2 - c^2 \xi^2}} \r e^{-\lambda t + t \sqrt{\lambda^2 - c^2 \xi^2}} + \l 1-\frac{\lambda}{\sqrt{\lambda^2 - c^2 \xi^2}} \r e^{-\lambda t - t \sqrt{\lambda^2 - c^2 \xi^2}} \right].
	\label{}
\end{align}
	In other words we have the following equality in distribution
	\begin{equation} 
		B \l c^2 \mathpzc{L}^{ \frac{1}{2} } (t) \r \,  \stackrel{\textrm{law}}{=} \, T \l \left| B(t) \right| \r, \qquad t>0.
		\label{}
	\end{equation}
\begin{proof}
First we show that the Fourier-Laplace transform of the law $w_{ \frac{1}{2} }^1 (x, t)$ of the process $W_1(t)  = S_1^2 \l c^2 \mathpzc{L}^{ \frac{1}{2}  } (t) \r = B \l c^2 \mathpzc{L}^{ \frac{1}{2}  } (t) \r $, $t>0$, is written as in \eqref{fourier laplace caso generale} for $\nu = \frac{1}{2}$, $\beta = 1$, $n=1$, and reads
\begin{equation} 
	\widehat{\widetilde{w_{ \frac{1}{2} }^1}} \l \xi, \mu \r \, = \, \frac{1+2\lambda \mu^{ -\frac{1}{2} }  }{\mu + 2\lambda \sqrt{\mu} + c^2 \xi^2}.
	\label{fl un mezzo}
\end{equation}
We have that
\begin{align}
	\widetilde{w^1_{ \frac{1}{2} }} (x, \mu) \, & = \, \int_0^\infty e^{-\mu t} \l \int_0^t p_{B} \l x, c^2s \r \, \mathpzc{l}_{ \frac{1}{2} } (s, t) \, ds \r \, dt \notag \\
	& = \, \int_0^\infty p_{B}(x, c^2 s) \, ds \int_s^\infty e^{-\mu t} \mathpzc{l}_{ \frac{1}{2} } (s, t) \, dt \notag \\
	& = \, \int_0^\infty p_{B} (x, c^2 s) \, ds \left[ \int_s^\infty e^{-\mu t}  \l \frac{2\lambda s \, e^{ -\frac{\l 2\lambda s \r^2}{4(t-s)} }}{\sqrt{4\pi \l t-s \r^3}} + 2\lambda \frac{e^{ -\frac{\l 2\lambda s \r^2}{4(t-s)} }}{\sqrt{\pi (t-s)}} \r dt \right] \notag \\
	& = \, \int_0^\infty p_{B} \l x, c^2 s \r \,  \l e^{-s \l \mu + 2\lambda \sqrt{\mu} \r} + 2\lambda \sqrt{\mu} e^{-s \l \mu + 2\lambda \sqrt{\mu} \r} \r \, ds \notag \\
	& = \, \int_0^\infty \frac{e^{-\frac{x^2}{4c^2s}}}{\sqrt{4\pi c^2s}}  e^{-s \l \mu + 2\lambda \sqrt{\mu} \r} \, ds + 2\lambda \mu^{ -\frac{1}{2} } \int_0^\infty \frac{e^{-\frac{x^2}{4c^2s}}}{\sqrt{4\pi c^2s}}  e^{-s \l \mu + 2\lambda \sqrt{\mu} \r} \, ds,
	\label{}
\end{align}
and thus taking the Fourier transform we get
\begin{align}
	\widehat{\widetilde{w^1_{ \frac{1}{2} }}} \l \xi, \mu \r \, & = \, \int_0^\infty e^{-sc^2 \xi^2} e^{-s \l \mu + 2\lambda \sqrt{\mu} \r} \, ds + 2\lambda \mu^{ -\frac{1}{2} } \int_0^\infty e^{-sc^2 \xi^2} e^{-s \l \mu + 2\lambda \sqrt{\mu} \r} \, ds \notag \\
	& = \, \frac{1+2\lambda \mu^{ -\frac{1}{2} }  }{\mu + 2\lambda \sqrt{\mu} + c^2 \xi^2}.
		\label{}
\end{align}

Now we are going to prove that the law $\mathpzc{w} (x, t)$ of the process $\mathpzc{W} (t)$, $t>0$, has Fourier-Laplace transform that coincides with \eqref{fl un mezzo}.
We have that
\begin{align}
	\mathpzc{w}(x, t) \, = \, \int_0^\infty p_T (x, s) \, p_{|B|} (s, t) \, ds,
	\label{}
\end{align}
and thus the Fourier transform of $\mathpzc{w} (x, t)$ reads
\begin{align}
	  \widehat{\mathpzc{w}} \l \xi, t \r \, =  \, & \int_{-\infty}^\infty e^{i\xi x} dx \int_0^\infty p_T \l x, s \r \, p_{|B|} (s, t) \, ds \,  \notag \\
	= \, & \frac{1}{2}  \int_0^\infty \left[  \l 1+\frac{\lambda}{\sqrt{\lambda^2 - c^2 \xi^2}} \r e^{-\lambda s + s \sqrt{\lambda^2 - c^2 \xi^2}} \right.  \notag \\
	 & \left. + \l 1-\frac{\lambda}{\sqrt{\lambda^2 - c^2 \xi^2}} \r e^{-\lambda s - s \sqrt{\lambda^2 - c^2 \xi^2}} \right] \, p_{|B|}(s, t) \, ds.
	\label{}
\end{align}
Passing now to the Laplace transform we have
\begin{align}
	\widetilde{\widehat{\mathpzc{w}}} (\xi, \mu) \,  = \, & \frac{1}{2} \int_0^\infty e^{-\mu t} \, dt   \int_0^\infty \left[  \l 1+\frac{\lambda}{\sqrt{\lambda^2 - c^2 \xi^2}} \r e^{-\lambda s + s \sqrt{\lambda^2 - c^2 \xi^2}} \;  \right. \notag \\
	&  + \left. \l 1-\frac{\lambda}{\sqrt{\lambda^2 - c^2 \xi^2}} \r e^{-\lambda s - s \sqrt{\lambda^2 - c^2 \xi^2}} \right] \, \frac{e^{-\frac{s^2}{4t}}}{\sqrt{\pi t}} \, ds \notag \\
	 = \, & \frac{1}{2}  \int_0^\infty \left[ \l 1+\frac{\lambda}{\sqrt{\lambda^2 - c^2 \xi^2}} \r e^{-\lambda s + s \sqrt{\lambda^2 - c^2 \xi^2}} \;  \right. \notag \\
	& + \left. \l 1-\frac{\lambda}{\sqrt{\lambda^2 - c^2 \xi^2}} \r e^{-\lambda s - s \sqrt{\lambda^2 - c^2 \xi^2}} \right] \, \frac{e^{ -s \sqrt{\mu} }}{\sqrt{\mu}}  \, ds \notag \\
	 = \, & \frac{1}{2\sqrt{\mu}} \left[ \l 1+\frac{\lambda}{\sqrt{\lambda^2 - c^2 \xi^2}} \r \l \frac{1}{\lambda + \sqrt{\mu} - \sqrt{\lambda^2 - c^2 \xi^2}} \r \;  \right. \notag \\
	&  + \left.  \l 1-\frac{\lambda}{\sqrt{\lambda^2 - c^2 \xi^2}} \r \l \frac{1}{\lambda + \sqrt{\mu} + \sqrt{\lambda^2 - c^2 \xi^2} } \r \right] \notag \\
	 = \, & \frac{\l \lambda + \sqrt{ \lambda^2 - c^2 \xi^2 } \r \l \lambda + \sqrt{\mu} + \sqrt{\lambda^2 - c^2 \xi^2} \r }{ \l 2\sqrt{\mu} \sqrt{\lambda^2 - c^2 \xi^2} \r \l \mu + 2\lambda \sqrt{\mu} + c^2 \xi^2 \r}  \notag \\
	&  + \frac{ \l \sqrt{\lambda^2 - c^2 \xi^2} - \lambda \r \l \lambda +  \sqrt{\mu} - \sqrt{\lambda^2 - c^2 \xi^2} \r }{ \l 2\sqrt{\mu} \sqrt{\lambda^2 - c^2 \xi^2} \r \l \mu + 2\lambda \sqrt{\mu} + c^2 \xi^2 \r} \notag \\
	= \, & \frac{1+ 2\lambda \mu^{ -\frac{1}{2} } }{\mu + 2\lambda \sqrt{\mu} + c^2 \xi^2},
	\label{}
\end{align}
which coincides with \eqref{fl un mezzo}.
\end{proof}
\end{te}
This shows that for each $t$ we have the following equality in distribution
\begin{equation} 
	T \l \lll B(t) \rrr \r \, \stackrel{\textrm{law}}{=} \, B \l c^2 \mathpzc{L}^{\frac{1}{2}} (t) \r, \qquad t>0,
	\label{role}
\end{equation}
where the role of the Brownian motion is interchanged in the two members of \eqref{role}. Thus, by suitably slowing down the time in \eqref{role}, we obtain the same distributional effect of a telegraph process taken at a Brownian time.
\begin{os} \normalfont
The probability distribution of the process
\begin{equation} 
	W_1(t) \, = \, B \l c^2 \mathpzc{L}^{ \frac{1}{2} } (t) \r, \qquad t>0,
	\label{}
\end{equation}
can be written as
\begin{align}
	w^1_{ \frac{1}{2} } (x, t) \, & = \, \frac{\lambda}{ c \pi} \int_0^t \frac{1}{\sqrt{s(t-s)}} e^{-\frac{x^2}{4c^2 s}-\frac{\lambda^2 s^2}{t-s}} \left[ \frac{s}{2(t-s)} +1 \right] \, ds \notag \\
	& = \, \frac{\lambda}{ c \pi} \int_0^t \frac{1}{\sqrt{s(t-s)}} e^{-\frac{x^2}{4c^2 s}-\frac{\lambda^2 s^2}{t-s}} \left[ \frac{1}{2} \l 1+\frac{t}{t-s} \r   \right] \, ds \notag \\
	& \stackrel{y=\lambda s}{=} \, \frac{\sqrt{\lambda}}{c\pi} \int_0^{\lambda t} e^{-\frac{\lambda x^2}{4c^2y}} \, e^{-\frac{y^2}{t-\frac{y}{\lambda}}} \frac{1}{\sqrt{y} \sqrt{t-\frac{y}{\lambda}}} \left[ \frac{1}{2} \l 1+\frac{t}{t-\frac{y}{\lambda}} \r \right] \, dy.
	\label{legge di B di l un mezzo}
\end{align}
Taking the limit for $c \to \infty$, $\lambda \to \infty$, $\frac{c^2}{\lambda} \to 1$, formula \eqref{legge di B di l un mezzo} becomes
\begin{equation} 
	\lim_{\begin{subarray}{c} \lambda, c \to \infty \\  \frac{c^2}{\lambda} \to 1  \end{subarray} } y_{ \frac{1}{2} }^1 \l x, t \r \, = \, 2 \int_0^\infty \frac{e^{-\frac{x^2}{4y}}}{\sqrt{4\pi y}} \frac{e^{-\frac{y^2}{t}}}{\sqrt{\pi t}} \, dy
	\label{}
\end{equation}
which coincides with the distribution of an iterated Brownian motion $B_1 \l \left| B_2(t) \right| \r$, $t>0$, with $B_j, j=1, 2$, independent Brownian motions.
From \eqref{legge di B di l un mezzo} we can see that the distribution of $W_1(t)$, $t>0$, has a bell-shaped structure.

Finally we show that the density $w_{ \frac{1}{2} }^1 \l x, t \r$ integrates to unity in force of the calculation
\begin{align}
	\int_{-\infty}^\infty w_{ \frac{1}{2} }^1 \l x, t \r \, dx \, & = \, \int_{-\infty}^\infty dx \int_0^t ds \, \frac{e^{-\frac{x^2}{4s}}}{\sqrt{4\pi s}} \, \mathpzc{l}_{ \frac{1}{2} } \l s, t \r \, = \,  \int_0^t ds \, \l \frac{\partial}{\partial s} \int_{ \frac{t-s}{\l 2\lambda \r^2} }^\infty \frac{s \, e^{-\frac{s^2}{4z}}}{\sqrt{4\pi z^3}} \, dz \r \notag \\
	& = \, \left[ \int_{ \frac{t-s}{\l 2\lambda \r^2} }^\infty \frac{s \, e^{-\frac{s^2}{4z}}}{\sqrt{4\pi z^3}} \, dz \right]_{s=0}^{s=t} \, = \, \int_0^\infty \frac{te^{-\frac{t^2}{4z}}}{\sqrt{4\pi z^3}} \, dz \, = \, 1.
	\label{}
\end{align}
In the intermediate step, formula \eqref{l un mezzo} has been applied.
\end{os}
\begin{os} \normalfont
The characteristic function of the process $T^{2\beta} (t)$, $t>0$, whose distribution satisfies
\begin{align}
\begin{cases}
	 \l \frac{\partial^2}{\partial t^2} + 2\lambda \frac{\partial}{\partial t} \r \, p_T^{2\beta}(x, t) \, = \, c^2 \frac{\partial^{2\beta}}{\partial |x|^{2\beta}} p_T^{2\beta}(x, t), \qquad 0 < \beta < 1, \beta \neq \frac{1}{2} \\
	 p_T^{2\beta}(x, 0) \, = \, \delta (x), \\
	\left. \frac{\partial}{\partial t} p_T^{2\beta}(x, t) \right|_{t=0} \, = \, 0,
	\end{cases}
	\label{}
\end{align}
reads
\small{
\begin{align}
	  & \mathbb{E}e^{i\xi T^{2\beta} (t)} \, = \notag \\
	  =  & \frac{e^{-\lambda t}}{2} \left[ \l 1+ \frac{\lambda}{\sqrt{\lambda^2 - c^2 |\xi |^{2\beta}}} \r  e^{t \sqrt{\lambda^2 - c^2 |\xi |^{2\beta}}}   + \l 1-\frac{\lambda}{\sqrt{\lambda^2 - c^2 |\xi |^{2\beta}}} \r e^{-t \sqrt{\lambda^2 - c^2 |\xi |^{2\beta}} } \right]
	\label{}
\end{align}
}
see \citet{zhao}. Therefore by performing the same steps as in theorem \eqref{teorema un mezzo} we prove that
\begin{equation} 
	S_1^{2\beta} \l \mathpzc{L}^{ \frac{1}{2} } (t) \r \, \stackrel{\textrm{law}}{=} \, T^{2\beta} \l \left| B(t) \right| \r, \qquad t>0.
	\label{}
\end{equation}
\end{os}

\subsection{The case $\nu = \frac{1}{3}$, convolutions of Airy functions}
We first recall that the totally positively-skewed stable process $H^{ \frac{1}{3} } (t)$, $t>0$ has law
\begin{equation} 
	\Pr \ll H^{ \frac{1}{3} } (t) \in dx \rr \, = \, \frac{t}{x \sqrt[3]{3x}} \textrm{Ai} \l \frac{t}{\sqrt[3]{3x}} \r \, dx, \qquad x >0, t>0,
	\label{legge di h un terzo}
\end{equation}
where $\textrm{Ai} (\cdot)$ is the Airy function.
Result \eqref{legge di h un terzo} can be obtained from the general series expansion of the stable law of order $\frac{1}{3}$ (see \citet{orsann} page 245) which reads
\begin{align}
	h_{ \frac{1}{3} } (x, 1) \, & = \, \frac{1}{3\pi} \sum_{k=0}^\infty (-1)^k \frac{\Gamma \l \frac{k+1}{3} \r}{k!} x^{-\frac{1}{3}(k+1)-1} \sin \l \frac{\pi}{3}(k+1) \r \notag \\
	& = \, \frac{1}{3\pi} \sum_{k=0}^\infty (-1)^k \frac{\Gamma \l \frac{k+1}{3} \r}{k!} x^{-\frac{k+1}{3}-1} (-1)^k \sin \l \frac{2\pi (k+1)}{3} \r \notag \\
	& = \, \frac{1}{3} \frac{3^{ \frac{2}{3} }}{x\sqrt[3]{x}} \textrm{Ai} \l \frac{1}{\sqrt[3]{3x}} \r \, = \, \frac{1}{x\sqrt[3]{3x}} \textrm{Ai} \l \frac{1}{\sqrt[3]{3x}} \r,
	\label{}
\end{align}
where we used formula (4.10) of \cite{orsann}, which reads
\begin{equation} 
	\textrm{Ai}(w) \, = \, \frac{3^{-\frac{2}{3}}}{\pi} \sum_{k=0}^\infty \l 3^{ \frac{1}{3} } w \r^k \frac{\sin \l \frac{2\pi (k+1)}{3} \r}{k!} \Gamma \l \frac{k+1}{3} \r.
	\label{}
\end{equation}
Since
\begin{equation} 
	H^{ \frac{1}{3} } (t) \, \stackrel{\textrm{law}}{=} \, t^3 H^{ \frac{1}{3} } (1),
	\label{}
\end{equation}
we have result \eqref{legge di h un terzo}. From the relatioship between $H^{ \frac{1}{3} } (t)$, $t>0$, and the inverse process $L^{ \frac{1}{3} } (t)$, $t>0$,
\begin{equation} 
	\Pr  \ll H^{ \frac{1}{3} } (t) <x \rr \, = \, \Pr \ll L^{ \frac{1}{3} } (x) >t \rr
	\label{}
\end{equation}
we extract the density of $L^{\frac{1}{3}} (x)$, $x>0$,
\begin{align}
	\frac{\Pr \ll L^{ \frac{1}{3} } (x) \in dt \rr}{dt}  \, & = \, -\frac{\partial}{\partial t} \int_0^x \frac{t}{s} \frac{1}{\sqrt[3]{3s}} \textrm{Ai} \l \frac{t}{\sqrt[3]{3s}} \r \, ds \notag \\
	& = \, -\int_0^x \frac{1}{s\sqrt[3]{3s}} \textrm{Ai} \l \frac{t}{\sqrt[3]{3s}} \r \, ds - \int_0^x \frac{t}{s\sqrt[3]{3s}} \textrm{Ai}^\prime \l \frac{t}{\sqrt[3]{3s}} \r \frac{ds}{\sqrt[3]{3s}}.
	\label{}
\end{align}
Since
\begin{equation} 
	\frac{\partial}{\partial s} \textrm{Ai} \l \frac{t}{\sqrt[3]{3s}} \r \, = \, - \frac{t}{3s \sqrt[3]{3s}} \textrm{Ai}^\prime \l \frac{t}{\sqrt[3]{3s}} \r
	\label{}
\end{equation}
we conclude that, for $x>0$, $t>0$,
\begin{align}
	l_{ \frac{1}{3} } (t, x) \,  = \, & \frac{\Pr \ll L^{ \frac{1}{3} } (x) \in dt \rr}{dt}   \notag \\
	 = \, & \int_0^x \frac{-1}{s\sqrt[3]{3s}} \textrm{Ai} \l \frac{t}{\sqrt[3]{3s}} \r \, ds + \int_0^x \frac{3}{\sqrt[3]{3s}} \frac{\partial}{\partial s} \textrm{Ai} \l \frac{t}{\sqrt[3]{3s}} \r \, ds \notag \\
	 = \, & \int_0^x \frac{-1}{s\sqrt[3]{3s}} \textrm{Ai} \l \frac{t}{\sqrt[3]{3s}} \r \, ds + \left[ \frac{3}{\sqrt[3]{3s}} \textrm{Ai} \l \frac{t}{\sqrt[3]{3s}} \r \right]_{s=0}^{s=x}   + \int_0^x \frac{ds}{s\sqrt[3]{3s}} \textrm{Ai} \l \frac{t}{\sqrt[3]{3s}} \r \,  \notag \\
	  = \, & \frac{3}{\sqrt[3]{3x}} \textrm{Ai} \l \frac{t}{\sqrt[3]{3x}} \r.
	\label{legge di elle un terzo}
\end{align}
In the last step we took into account the asymptotic expansion 7.2.19 of \citet{bleis}.

With similar calculation we obtain the law $h_{ \frac{2}{3} } (x, t)$ of the process $H^{\frac{2}{3}} (t)$, $t>0$, which is expressed in terms of Airy function. From the general series expression of the stable law (see \cite{orsann}) we have that,
\begin{align}
	& \qquad h_{ \frac{2}{3} } (x, 1) \,  = \notag \\
	& = \, \frac{2}{3\pi} \sum_{k=0}^\infty (-1)^k \frac{\Gamma \l \frac{2}{3} (k+1) \r}{k!} x^{ -\frac{2}{3} (k+1) -1 } \sin \l \frac{2\pi}{3}(k+1) \r \notag \\
	& = \, \frac{2}{3\pi \sqrt{\pi}} \sum_{k=0}^\infty \frac{(-1)^k}{k!} \frac{x^{ -\frac{2}{3} (k+1)-1 }}{2^{1-\frac{2}{3}(k+1)}} \Gamma \l \frac{k+1}{3} \r \sin \l \frac{2\pi}{3} (k+1) \r \int_0^\infty dw \, e^{-w} w^{\frac{k+1}{3}+\frac{1}{2}-1} \notag \\
	& = \, \frac{1}{x} \sqrt[3]{\frac{2^2}{3x^2}} \frac{1}{\sqrt{\pi}} \int_0^\infty e^{-w} w^{ -\frac{1}{6} } \, \textrm{Ai} \l -\sqrt[3]{\frac{2^2w}{3x^2}} \r \, dw,
	\label{}
\end{align}
and thus, in force of the fact that $H^{ \frac{2}{3} } (t) \stackrel{\textrm{law}}{=} t^{ \frac{3}{2} } H^{\frac{2}{3}} (1)$,
\begin{equation} 
h_{ \frac{2}{3} } (x, t) \, = \, \frac{t}{\sqrt{\pi}} \frac{1}{x} \int_0^\infty dw \, e^{-w} w^{-\frac{1}{6}} \, \sqrt[3]{\frac{2^2}{3x^2}} \, \textrm{Ai} \l -t\sqrt[3]{\frac{2^2w}{3x^2}} \r.
	\label{legge di h due terzi}
\end{equation}
\begin{os} \normalfont
\label{remarkintegr1airy}
We check that the distribution \eqref{legge di h due terzi} integrates to unity. We have that
\begin{align}
	& \qquad \int_0^\infty h_{\frac{2}{3}} (x, t) \, dx \,  = \notag \\
	& = \, \frac{t}{\sqrt{\pi}}  \int_0^\infty dw \, e^{-w}  w^{-\frac{1}{6}} \sqrt[3]{\frac{2^2}{3}} \int_0^\infty dx \, x^{-\frac{2}{3}-1} \, \textrm{Ai} \l -t\sqrt[3]{\frac{2^2w}{3x^2}} \r \notag \\
	& \stackrel{y = x^{-\frac{2}{3}} t \sqrt[3]{\frac{2^2w}{3}} }{=} \,  \frac{t}{\sqrt{\pi}}  \int_0^\infty dw \, e^{-w}  w^{-\frac{1}{6}} \sqrt[3]{\frac{2^2}{3}} \frac{3}{2} \l t \sqrt[3]{\frac{2^2w}{3}} \r^{-1} \int_0^\infty dy \, \textrm{Ai} \l -y \r \notag \\
	& = \, \frac{1}{\sqrt{\pi}}  \int_0^\infty dw \, e^{-w}  w^{-\frac{1}{6}} \sqrt[3]{\frac{2^2}{3}}  \l  \sqrt[3]{\frac{2^2w}{3}} \r^{-1} \notag \\
	& = \, \frac{1}{\sqrt{\pi}}  \int_0^\infty dw \, e^{-w}  w^{-\frac{1}{6}-\frac{1}{3}} \, = \, \frac{1}{\sqrt{\pi}}  \int_0^\infty dw \, e^{-w}  w^{\frac{1}{2}-1} \, = \, 1,
	\label{}
\end{align}
where we used the fact that
\begin{equation} 
	\int_0^\infty dy \, \textrm{Ai} (-y) \,  = \, \frac{2}{3}.
	\label{}
\end{equation}
\end{os}
For the law of the process $L^{\frac{2}{3}} (x)$, $x>0$, we therefore have that
\begin{align}
	\Pr \ll L^{ \frac{2}{3} } (x) <t \rr \, & = \, \Pr \ll H^{ \frac{2}{3} } (t) > x \rr \notag \\
	& = \, \int_0^\infty \int_x^\infty \frac{t}{\sqrt{\pi}} \frac{1}{z} \sqrt[3]{\frac{2^2}{3z^2}} \, \textrm{Ai} \l -t \sqrt[3]{\frac{2^2w}{3z^2}} \r \, e^{-w} w^{ -\frac{1}{6} } \, dw \, dz 
	\label{}
\end{align}
and thus
\begin{align}
	 l_{ \frac{2}{3} } \l t, x \r \, = \, & \int_0^\infty \int_x^\infty \frac{dw \, dz}{z \sqrt{\pi}} \sqrt[3]{\frac{2^2}{3z^2}} \, \textrm{Ai} \l -t \sqrt[3]{\frac{2^2 w}{3z^2}} \r \, e^{-w} w^{-\frac{1}{6}} \notag \\
	 &  - \int_0^\infty \int_x^\infty \frac{t}{z \sqrt{\pi}} \sqrt[3]{\frac{2^2}{3z^2}} \sqrt[3]{\frac{2^2w}{3z^2}} \, \textrm{Ai}^\prime \l -t \sqrt[3]{\frac{2^2w}{3z^2}} \r \, dz \, dw \notag \\
	 = \, & \int_0^\infty \int_x^\infty \frac{dw \, dz}{z \sqrt{\pi}} \sqrt[3]{\frac{2^2}{3z^2}} \, \textrm{Ai} \l -t \sqrt[3]{\frac{2^2 w}{3z^2}} \r \, e^{-w} w^{-\frac{1}{6}} \,  \notag \\
	& - \frac{3}{2} \int_x^\infty \int_0^\infty \frac{1}{\sqrt{\pi}} \sqrt[3]{\frac{2^2}{3z^2}} e^{-w} w^{-\frac{1}{6}} \frac{\partial}{\partial z} \textrm{Ai} \l -t \sqrt[3]{\frac{2^2w}{3z^2}} \r \, dw \, dz \notag \\
	 = \, & \int_0^\infty \int_x^\infty \frac{dw \, dz}{z \sqrt{\pi}} \sqrt[3]{\frac{2^2}{3z^2}} \, \textrm{Ai} \l -t \sqrt[3]{\frac{2^2 w}{3z^2}} \r \, e^{-w} w^{-\frac{1}{6}}  \notag  \\  
	&  - \left[ \frac{3}{2 \sqrt{\pi}} \int_0^\infty dw  \sqrt[3]{\frac{2^2}{3z^2}} e^{-w} w^{-\frac{1}{6}} \, \textrm{Ai} \l -t\sqrt[3]{\frac{2^2w}{3z^2}} \r \right]_{z=x}^{z=\infty}  \notag \\
	&  - \int_0^\infty \int_x^\infty \frac{dw \, dz}{z \sqrt{\pi}} \sqrt[3]{\frac{2^2}{3z^2}} \, \textrm{Ai} \l -t \sqrt[3]{\frac{2^2 w}{3z^2}} \r \, e^{-w} w^{-\frac{1}{6}} \notag \\
	 = \, & \frac{3}{2\sqrt{\pi}} \int_0^\infty  \sqrt[3]{\frac{2^2}{3x^2}} e^{-w} w^{-\frac{1}{6}} \, \textrm{Ai} \l -t \sqrt[3]{\frac{2^2w}{3x^2}} \r \, dw.
	\label{legge di l due terzi}
\end{align}
For checking that \eqref{legge di l due terzi} integrates to unity one can perform calculation similar to that of Remark \ref{remarkintegr1airy}.

Now we have all the information to get the distribution of the process $\mathpzc{L}^{ \frac{1}{3} } (t)$, $t>0$, by means of formula \eqref{legge esplicita di lbold per ogni nu}.
We have that
\begin{align}
	\mathpzc{l}_{ \frac{1}{3} } (x, t) \,  = \, & \frac{\Pr \ll \mathpzc{L}^{ \frac{1}{3} } (t) \in dx \rr}{dx} \, \notag \\
	 = \, & \int_0^t l_{\frac{2}{3}} (x, t-s) \, h_{\frac{1}{3}} (s, 2\lambda x) ds + 2\lambda \int_0^t l_{\frac{1}{3}} (2\lambda x, s) \, h_{\frac{2}{3}} (t-s, x) \, ds \notag \\
	 = \, & \int_0^t ds \, \left[ \frac{3}{2\sqrt{\pi}} \int_0^\infty dw \sqrt[3]{\frac{2^2}{3(t-s)^2}} e^{-w} w^{-\frac{1}{6}} \, \textrm{Ai} \l -x \sqrt[3]{\frac{2^2w}{3(t-s)^2}} \r \, dw  \right] \bm{\cdot} \notag \\
	& \bm{\cdot} \frac{2\lambda x}{s \sqrt[3]{3s}} \, \textrm{Ai} \l \frac{2\lambda x}{\sqrt[3]{3s}} \r  + 2\lambda \int_0^t ds \frac{3}{\sqrt[3]{3s}} \, \textrm{Ai} \l \frac{2\lambda x}{\sqrt[3]{3s}} \r \bm{\cdot} \notag \\
	& \bm{\cdot} \frac{s}{\sqrt{\pi} (t-s)} \int_0^\infty dw \, e^{-w} w^{-\frac{1}{6}} \sqrt[3]{\frac{2^2}{3(t-s)^2}} \, \textrm{Ai} \l -x \sqrt[3]{\frac{2^2w}{3(t-s)^2}} \r \notag \\
	 = \, & \frac{2\lambda}{\sqrt{\pi}} \int_0^t ds \int_0^\infty dw \, e^{-w} w^{-\frac{1}{6}} \, \textrm{Ai} \l -x \sqrt[3]{\frac{2^2w}{3(t-s)^2}} \r \, \textrm{Ai} \l \frac{2\lambda x}{\sqrt[3]{3s}} \r \bm{\cdot} \notag \\
	&  \bm{\cdot} \frac{3}{\sqrt[3]{3s}} \sqrt[3]{\frac{2^2}{3(t-s)^2}} \left[ \frac{x}{2s} + \frac{s}{t-s} \right].
	\label{4633}
\end{align}
Result \eqref{4633} permits us to write explicitly the solution of the fractional telegraph equation \eqref{1.1} for $\nu = \frac{1}{3}$, $\beta = 1$ and $n = 1$, as
\begin{equation} 
	w_{\frac{1}{3}}^1 (x, t) \, = \, \int_0^\infty \frac{e^{-\frac{x^2}{4c^2 s}}}{\sqrt{4 \pi c^2  s}} \, \mathpzc{l}_{\frac{1}{3}} (s, t) \, ds, \qquad x \in \mathbb{R}, t>0.
	\label{}
\end{equation}

\subsection{The planar case}
Let us consider the planar process 
\begin{equation}
	\bm{T}(t) \, = \, \l X(t), Y(t) \r, \qquad t>0,
\label{}
\end{equation}
with infinite directions and finite velocity $c$, investigated in \citet{orsdegr}, which has probability law (see formula 1.2 therein)
     	\begin{align}
   	   	r(x,y,t)= \frac{\lambda}{2\pi c} \frac{e^{-\lambda t +\frac{\lambda}{c}\sqrt{c^2 t^2 - (x^2 + y^2)}}}{\sqrt{c^2 t^2 - (x^2+ y^2)}}, \qquad x^2+y^2 < c^2t^2, t>0,
     	\label{eq:2}
     	\end{align}
which satisfies the telegraph equation
\begin{align}
    	\l \frac{\partial^2 }{\partial t^2}+2\lambda \frac{\partial }{\partial t} \r r(x,y,t) \, = \, c^2\left(\frac{\partial^2}{\partial x^2}+ \frac{\partial^2}{\partial y^2}\right)\, r(x,y,t).
       	\label{eq:6}
   	\end{align}
The distribution of $\bm{T} (t)$, $t>0$, has a singular component uniformly distributed on the circle $\partial C_{ct} = \ll (x, y) \in \mathbb{R}^2 : x^2+y^2 = c^2t^2 \rr$ with probability mass equal to $e^{-\lambda t}$.
The process $\bm{T} (t)$, $t>0$, describes a random motion where directions change at Poisson paced times and the orientation of each segment of the sample paths is uniform in $[0, 2\pi)$.

     Let $q(x,y,t)$ be the distribution obtained by means of the composition of the process $\bm{T}(t)$ with a reflecting Brownian motion with law
	\begin{align}
     	p_{|B|}(s,t)\, = \, \frac{e^{\frac{-s^2}{4t}}}{\sqrt{\pi t}}, \qquad t>0, s>0,
     	\label{eq:3}
     	\end{align}
which satisfies the equation
       	\begin{align}
    	\frac{^C\partial^{\frac{1}{2}}  }{\partial t^{\frac{1}{2}}} p_{|B|}(s,t) \, = \, - \frac{\partial }{\partial s} p_{|B|}(s,t)
    	\label{eq:4}
    	\end{align}
     and also
     	\begin{align}
     	\frac{\partial }{\partial t}  p_{|B|}(s,t) \, = \,  \frac{\partial^2 }{\partial s^2} p_{|B|}(s,t)
     	\label{eq:5}
     	\end{align}
     	We have the following theorem.
     	\begin{te}
     The law of the composition
     \begin{equation} 
     	\bm{Q} (t) \, = \, \bm{T} \l \left| B (t) \right| \r, \qquad t>0
     	\label{}
     \end{equation}
    written as
     	\begin{align}
     	q(x,y,t)= \int_0^\infty \, r(x,y,s) \, p_{|B|}(s,t) \, ds,
     	\label{eq:1}
     	\end{align}
     	satisfies the $2$-dimensional time-fractional equation
     	\begin{equation} 
     		\l \frac{\partial}{\partial t}  + 2\lambda \frac{^C\partial^{\frac{1}{2}}}{\partial t^{\frac{1}{2}}} \r q \l x, y, t \r \, = \, c^2 \l \frac{\partial^2}{\partial x^2} + \frac{\partial^2}{\partial y^2} \r q \l x, y, t \r, \qquad x, y \in \mathbb{R}, t>0,
     		\label{equazione caso planare}
     	\end{equation}
     	subject to the initial condition $q(x, y, 0) \, = \, \delta (x, y)$.

   \begin{proof}
By considering \eqref{eq:1} and \eqref{eq:4} we can write
   	\begin{align}
    	\frac{^C\partial^\frac{1}{2} }{\partial t^\frac{1}{2}} q(x,y,t) \, & = \, \int_0^\infty \, r(x,y,s)\,   \frac{^C\partial^\frac{1}{2} }{\partial t^\frac{1}{2}} p_{|B|}(s,t) \, ds  \notag \\
    	& = \, \int_0^\infty \, r(x,y,s) \left( - \frac{\partial }{\partial s} p_{|B|}(s,t) \right)\, ds  \notag \\
    	& = \, \left[-p_{|B|}(s,t)\, r(x,y,s)\right]_{s=0}^{s=\infty} + \int_0^\infty \, p_{|B|}(s,t) \frac{\partial }{\partial s} r(x,y,s) \, ds.    	      	  \label{eq:7}		    	
    	\end{align}
In the previous step it must be taken into account that the boundary $\partial C_{cs}$ is excluded. From \eqref{eq:1} and \eqref{eq:5} we have that
    	\begin{align}    	
    	\frac{\partial }{\partial t} q(x,y,t) \, & = \, \int_0^\infty \, r(x,y,s) \,  \frac{\partial }{\partial t} p_{|B|}(s,t)  \, ds  \, = \,  \int_0^\infty \, r(x,y,s) \,  \frac{\partial^2 }{\partial s^2} p_{|B|}(s,t)  \, ds  \notag \\
    	& = \, \left[ r(x,y,s) \, \frac{\partial }{\partial s} p_{|B|}(s,t) \right]_{s=0}^{s=\infty} -  \int_0^\infty \, \frac{\partial }{\partial s} r(x,y,s) \, \frac{\partial }{\partial s} p_{|B|}(s,t) \, ds  \notag \\
    	& = \, -  \left[p_{|B|}(s,t) \, \frac{\partial }{\partial s} r(x,y,s) \right]_{s=0}^{s=\infty} +  \int_0^\infty \, p_{|B|}(s,t) \, \frac{\partial^2 }{\partial s^2} r(x,y,s)   \,ds.
    	\label{eq:8}
    	\end{align}
Thus, by looking at \eqref{eq:6}, \eqref{eq:7} and \eqref{eq:8} we obtain	    
 		\begin{align}
 		&  \frac{\partial }{\partial t} q(x,y,t) + 2\lambda \frac{^C\partial^\frac{1}{2} }{\partial t^\frac{1}{2}} q(x,y,t)  \, = \notag \\ = \, &\int_0^\infty p_{|B|} (s, t) \left[ \frac{\partial^2}{\partial s^2} r(x, y, s) + 2\lambda \frac{\partial}{\partial s} r(x, y, s) \right] \, ds \notag \\
 		 = \, & \int_0^\infty p_{|B|} (s, t) \, c^2 \l \frac{\partial^2}{\partial x^2} + \frac{\partial^2}{\partial y^2}  \r r(x, y, s) \, ds \, = \, c^2 \left( \frac{\partial^2}{\partial x^2}+\frac{\partial^2}{\partial y^2} \right) q(x,y,t).
 		\label{eq:9}
		\end{align} 		   
	which means that $q(x,y,t)$ satisfies equation \eqref{equazione caso planare}.
\end{proof}	
\end{te}

	It is easy to show that the process $\bm{Q}(t) = \bm{T} \l \left| B(t) \right| \r$, $t>0$, has not the same law of the process $\bm{W}_2 (t) = \bm{B}_2 \l c^2 \mathpzc{L}^{\frac{1}{2}} (t) \r$, $t>0$. However it is possible to construct a planar process, say $\bm{\mathfrak{T}} (t)$, $t>0$ (which is a slightly different version of $\bm{T}(t)$, $t>0$) composed with a suitable "time process" which has the same distribution as $\bm{W}_2 (t)$, $t>0$. The planar random motion $\bm{\mathfrak{T}} (t)$, $t>0$, with distribution
	\begin{equation} 
		\mathfrak{r} (x, y, t) \, = \, \frac{\lambda \, e^{-\lambda t}}{2 \pi c} \left[ \frac{e^{\frac{\lambda}{c} \sqrt{c^2t^2 - \l x^2+y^2 \r}}+ e^{-\frac{\lambda}{c} \sqrt{c^2t^2-\l x^2+y^2 \r}}}{\sqrt{c^2t^2 - \l x^2+y^2 \r}}  \right],
		\label{478}
	\end{equation}
	where $(x, y) \in C_{ct} = \ll (x, y) : x^2+y^2 < c^2t^2 \rr$, can be constructed starting from the model dealt with in \citet{orsdegr}. The distribution is based on the solution to the planar telegraph equation
\begin{equation} 
	\l \frac{\partial^2}{\partial t^2} + 2\lambda \frac{\partial}{\partial t} \r \mathfrak{r} (x, y, t) \, = \, c^2 \l \frac{\partial^2}{\partial x^2} + \frac{\partial^2}{\partial y^2} \r \mathfrak{r} (x, y, t),
		\label{}
\end{equation}
namely
	\begin{equation} 
		\mathfrak{r} \l x, y, t \r \, = \,  \frac{e^{-\lambda t}}{\sqrt{c^2t^2 - \l x^2+y^2 \r}} \left[ Ae^{\frac{\lambda}{c} \sqrt{c^2 t^2 - \l x^2+y^2 \r} }+B e^{-\frac{\lambda}{c} \sqrt{c^2t^2- \l x^2+y^2 \r}} \right],
			\label{complete}
	\end{equation}
with $A = B = \frac{\lambda}{2\pi c}$ and thus we can easily check that
 \begin{equation} 
 	\iint_{C_{ct}} dx \, dy \, \mathfrak{r} \l x, y, t \r \, = \, 1-e^{-2\lambda t}.
 	\label{}
 \end{equation}
We take a particle starting from the origin, moving at finite velocity $c$, and changing direction (chosen with uniform distribution) at Poisson times and neglect displacements started off by even-labelled times. The sample paths of this motion are constructed by piecing together only odd-order displacements of the planar motion $\bm{T}(t)$, $t>0$. The process just described has distribution \eqref{478} as shown below
\begin{align}
	& \mathfrak{r} (x, y, t) = \notag \\
	 = \, & \frac{\Pr \ll \bm{\mathfrak{T}} (t)  \in d\bm{x}  \rr}{d\bm{x}} \, = \,  \frac{\lambda \, e^{-\lambda t}}{2 \pi c} \left[ \frac{e^{\frac{\lambda}{c} \sqrt{c^2t^2 - \l x^2+y^2 \r}}+ e^{-\frac{\lambda}{c} \sqrt{c^2t^2 - \l x^2+y^2 \r}}}{\sqrt{c^2t^2 - \l x^2+y^2 \r}} \right] \notag \\
	 = \, & \frac{\lambda^2}{c^2} \frac{1}{\pi} e^{-\lambda t} \left[ \sum_{k=0}^\infty \l \frac{\lambda}{c} \sqrt{c^2t^2 - \l x^2+y^2 \r} \r^{2k-1} \frac{1}{(2k)!} \right] \notag \\
	 = \, & \frac{\lambda^2}{c^2} \frac{1}{\pi} \sum_{k=0}^\infty \l \frac{\lambda}{c} \r^{2k-1} (2k+1) \l c^2t^2- \l x^2+y^2 \r \r^{k-\frac{1}{2}} \frac{e^{-\lambda t}}{(2k)! (2k+1)} \frac{(\lambda t)^{2k+1}}{(\lambda t)^{2k+1}} \notag \\
	 = \, & 2 \sum_{k=0}^\infty \Pr \ll X(t) \in dx, Y(t) \in dy | N(t) = 2k+1 \rr \, e^{-\lambda t} \frac{(\lambda t)^{2k+1}}{(2k+1)!} \notag \\
	 = \, & 2 \sum_{k=0}^\infty \Pr \ll \bm{T}(t) \in d\bm{x} | N(t) = 2k+1 \rr \, e^{-\lambda t} \frac{(\lambda t)^{2k+1}}{(2k+1)!},
	\label{475}
\end{align}
where, for $x^2+y^2 < c^2t^2$ (see \cite{orsdegr}),
\begin{equation}
	\frac{\Pr \ll X(t) \in dx, Y(t) \in dy | N(t) = n \rr}{dx \, dy} \, = \, \frac{n}{2n (ct)^n} \l c^2t^2 - \l x^2+y^2 \r \r^{\frac{n}{2}-1},
	\label{}
\end{equation}
and
\begin{equation} 
	2e^{-\lambda t} \sum_{k=0}^\infty \frac{(\lambda t)^{2k+1}}{(2k+1)!} \, = \, \sum_{k=0}^\infty 2 \Pr \ll N(t) \, = \, 2k+1 \rr \, = \, 1-e^{-2\lambda t}.
	\label{476}
\end{equation}
The factor $2$ appearing in \eqref{475} and \eqref{476} can be interpreted as follows. The displacements generated by an even number of Poisson events are disregarded and replaced by displacements produced by an odd number of deviations. Therefore, odd-order Poisson events ignite twice the displacements considered in \eqref{475}.
\begin{te}
The composition with distribution
\begin{equation} 
	\mathfrak{q} (x, y, t) \, = \, \int_0^\infty ds \, \mathfrak{r} \l x, y, s \r \left[ p_{|B|} \l s, t \r + \frac{1}{2\lambda} \frac{\partial^{\frac{1}{2}}}{\partial t^{\frac{1}{2}}} p_{|B|} \l s, t \r \right],
	\label{processo finito composto}
\end{equation}
which satisfies the time-fractional equation
\begin{equation} 
	\l \frac{\partial}{\partial t} + 2\lambda \frac{^C\partial^{\frac{1}{2}}}{\partial t^{\frac{1}{2}}} \r \mathfrak{q}(x, y, t) \, = \, c^2 \l \frac{\partial^2}{\partial x^2} + \frac{\partial^2}{\partial y^2} \r \mathfrak{q}(x, y, t),
	\label{timw frac bid}
\end{equation}
has the same law of the process $\bm{W}_2 (t) = \bm{B}_2 \l c^2 \mathpzc{L}^{\frac{1}{2}} (t) \r$.
\begin{proof}
We begin by evaluating the Fourier-Laplace transform of \eqref{processo finito composto}.
\begin{align}
	& \widehat{\widetilde{\mathfrak{q}}} (\xi, \alpha, \mu) \notag \\ 
	 = \, & \int_0^\infty ds \int_0^\infty dt \, e^{-\mu t} \int_{C_{ct}} dx \, dy \, e^{i\xi x + i \alpha y}  \mathfrak{r} (x, y, s) \left[p_{|B|} (s, t) + \frac{1}{2\lambda} \frac{\partial^{\frac{1}{2}}}{\partial t^{\frac{1}{2}}} p_{|B|} (s, t) \right] \notag \\
	 = \, & \frac{2\lambda  + \sqrt{\mu}}{2\lambda \sqrt{\mu}} \int_0^\infty ds \int_{C_{ct}} dx \, dy \, e^{i\xi x + i \alpha y} \, \mathfrak{r} \l x, y , s \r \, e^{-s\sqrt{\mu}}.
	\label{from questa}
\end{align}
Now we need the Fourier transform of the law $\mathfrak{r} (x, y, t)$ of the process $\bm{\mathfrak{T}} (t)$, $t>0$, which reads
\begin{align}
	&\widehat{\mathfrak{r}} (\xi, \alpha, t ) = \notag \\
	 = \, & \frac{\lambda \, e^{-\lambda t}}{2 \pi c} \iint_{C_{ct}} e^{i\xi x + i \alpha y} \left[ \frac{e^{\frac{\lambda}{c} \sqrt{c^2t^2 - \l x^2+y^2 \r} }+ e^{-\frac{\lambda}{c}\sqrt{c^2t^2 - \l x^2+y^2 \r}}}{\sqrt{c^2t^2-\l x^2+y^2 \r} } \right] dx \, dy \notag \\
	 = \, & \frac{\lambda \, e^{-\lambda t}}{2 \pi c} \int_0^{2\pi} d \theta \int_0^{ct} d\rho \, \rho e^{i\rho \l \xi \cos \theta + \alpha \sin \theta \r} \frac{\lambda}{c} \frac{e^{\frac{\lambda}{c} \sqrt{c^2t^2 - \rho^2}}+e^{-\frac{\lambda}{c}} \sqrt{c^2t^2 - \rho^2}}{\sqrt{c^2t^2-\rho^2}} \notag \\
	 = \, & \frac{2\lambda^2 e^{-\lambda t}}{c^2} \int_0^{ct} \rho \sum_{m=0}^\infty \l \frac{\lambda}{c} \sqrt{c^2t^2-\rho^2} \r^{2m-1} \frac{1}{(2m)!} \, J_0 \l \rho \sqrt{\xi^2 + \alpha^2} \r \, d\rho \notag \\
	 = \, & \frac{2\lambda e^{-\lambda t}}{c} \sum_{m=0}^\infty \frac{\l \frac{\lambda}{c} \r^{2m}}{(2m)!} \sum_{k=0}^\infty  \frac{(-1)^k \l \frac{\sqrt{\xi^2 + \alpha^2}}{2} \r^{2k}}{(k!)^2} \bm{\cdot} \int_0^{ct}  \l c^2t^2 - \rho^2 \r^{m-\frac{1}{2}} \rho^{2k+1} \, d\rho \notag \\
	 = \, & \frac{2\lambda e^{-\lambda t}}{c} \sum_{m=0}^\infty  \frac{\l \frac{\lambda}{c} \r^{2m}}{(2m)!} \sum_{k=0}^\infty    \frac{(-1)^k \l \frac{\sqrt{\xi^2 + \alpha^2}}{2} \r^{2k} }{2(k!)^2 \, (ct)^{-(2m+2k+1)}} \int_0^1 y^k \l 1-y \r^{m-\frac{1}{2}} dy \notag \\
	 = \, & \frac{\lambda}{c} e^{-\lambda t} \sum_{m=0}^\infty \l \frac{\lambda}{c} \r^{2m} \frac{1}{(2m)!} \sum_{k=0}^\infty (-1)^k \l \frac{\sqrt{\xi^2 + \alpha^2}}{2} \r^{2k} \frac{(ct)^{2m+2k+1} \Gamma \l m+\frac{1}{2} \r}{ k! \Gamma \l k+m+1+\frac{1}{2} \r}.
	\label{}
\end{align}
Thus, from \eqref{from questa}, we have that 
\begin{align}
	\widetilde{\widehat{\mathfrak{q}}} \l \xi, \alpha, \mu \r \, & = \,  \frac{2\lambda  + \sqrt{\mu}}{2\lambda \sqrt{\mu}} \int_0^\infty ds \, \widehat{\mathfrak{r}} (\xi, \alpha, t) \, e^{-s\sqrt{\mu}} \, = \frac{1+2\lambda \mu^{-\frac{1}{2}}}{\mu + 2\lambda \sqrt{\mu } + c^2 \l \xi^2 + \alpha^2 \r}
	\label{fourier laplace proc plan comp}
\end{align}
in force of the calculation
\begin{align}
 & \int_0^\infty ds \, \widehat{\mathfrak{r}} \l \xi, \alpha, s \r e^{-s\sqrt{\mu}} \,  =  \notag \\
 = \, & \frac{\lambda}{c} \int_0^\infty ds \, e^{-\lambda s} \sum_{m=0}^\infty  \frac{\lambda^{2m}}{c^{2m}(2m)!} \sum_{k=0}^\infty    \frac{(-1)^k \l \frac{\sqrt{\xi^2+\alpha^2}}{2} \r^{2k}\Gamma \l m+\frac{1}{2} \r}{k!\, (cs)^{-(2m+2k+1)} \, \Gamma \l k+m+1+\frac{1}{2} \r}  e^{-s\sqrt{\mu}} \notag \\
 = \, &  \lambda \sum_{m=0}^\infty \frac{\sqrt{\pi} 2^{1-2m} \Gamma \l 2m \r}{\lambda^{-2m}(2m)! \Gamma (m)} \sum_{k=0}^\infty  \frac{(-1)^k \l \frac{\sqrt{\xi^2 + \alpha^2}}{2} \r^{2k} c^{2k}}{k!\Gamma \l k+m+1+\frac{1}{2} \r} \int_0^\infty e^{-s \l \lambda + \sqrt{\mu} \r} \, s^{2m+2k+1}ds   \notag \\
 = \, &   \frac{\lambda}{2\l \lambda + \sqrt{\mu} \r^2} \sum_{m=0}^\infty \frac{\lambda^{2m} \sqrt{\pi} 2^{1-2m}}{m! \l \lambda + \sqrt{\mu} \r^{2m}} \sum_{k=0}^\infty  \frac{(-1)^k  \l \frac{\sqrt{\xi^2 + \alpha^2}}{2} \r^{2k} }{k! \l \lambda + \sqrt{\mu} \r^{2k}c^{-2k}} \frac{\Gamma \l 2k+2m+2 \r}{\Gamma \l k+m+1+\frac{1}{2} \r} \notag \\
 = \, & \frac{ \sqrt{\pi} \lambda}{2\l \lambda + \sqrt{\mu} \r^2} \sum_{m=0}^\infty \frac{\lambda^{2m} 2^{1-2m} }{m! \l \lambda + \sqrt{\mu} \r^{2m}} \sum_{k=0}^\infty \frac{(-1)^k  \l \frac{\sqrt{\xi^2+\alpha^2}}{2} \r^{2k} }{k! \l \lambda + \sqrt{\mu} \r^{2k} c^{-2k}} \frac{\Gamma \l k+m+1 \r}{2^{1-2(k+m+1)} \sqrt{\pi}} \notag \\
 = \, & \frac{2\lambda }{\l \lambda + \sqrt{\mu} \r^2} \sum_{m=0}^\infty \frac{\lambda^{2m}}{m! \l \lambda + \sqrt{\mu} \r^{2m}} \sum_{k=0}^\infty  \frac{(-1)^k\l \sqrt{\xi^2 + \alpha^2} \r^{2k}}{k! \l \lambda + \sqrt{\mu} \r^{2k}c^{-2k}}    \int_0^\infty e^{-u} u^{k+m} \, du \notag \\
 = \, & \frac{2\lambda }{\l \lambda + \sqrt{\mu} \r^2} \int_0^\infty du \, e^{u \frac{\lambda^2}{\l \lambda + \sqrt{\mu} \r^2}-u\frac{c^2 \l \xi^2 + \alpha^2 \r}{\l \lambda + \sqrt{\mu} \r^2}-u} \, = \, \frac{\frac{2\lambda }{\l \lambda + \sqrt{\mu} \r^2}}{1-\frac{\lambda^2}{\l \lambda + \sqrt{\mu} \r^2} + \frac{c^2 \l \xi^2 + \alpha^2 \r}{\l \lambda + \sqrt{\mu} \r^2}}  \notag \\
 = \, & \frac{2\lambda }{\l \lambda + \sqrt{\mu} \r^2- \lambda^2 + c^2 \l \xi^2 + \alpha^2 \r} \, = \, \frac{2\lambda }{\mu + 2\lambda \sqrt{\mu} + c^2 \l \xi^2 + \alpha^2 \r }.
\label{}
\end{align}

The Fourier-Laplace transform of the law of the process $\bm{B}_2 \l c^2 \mathpzc{L}^{\frac{1}{2}} (t) \r$ is written as in \eqref{fourier laplace caso generale} for $n=2$, $\beta = 1$ and $\nu = \frac{1}{2}$ as the following calculation shows
\begin{align}
\widehat{\widetilde{w_{\frac{1}{2}}^1}} \l \xi, \alpha, t \r \, & = \, \int_0^\infty \widehat{p_{\bm{B}}} \l \xi, \alpha, c^2 s \r \, \widetilde{\mathpzc{l}_{\frac{1}{2}}} \l s, \mu \r \, ds \notag \\
& = \l 1 + 2\lambda \mu^{-\frac{1}{2}} \r \int_0^\infty e^{-\mu s - \l \xi^2 + \alpha^2 \r c^2 s} \left[ e^{-2\lambda s \sqrt{\mu}} + 2\lambda \frac{e^{-2\lambda s \sqrt{\mu}}}{\sqrt{\mu}} \right] \, ds \notag \\
& = \, \frac{1+2\lambda \mu^{-\frac{1}{2}}}{2\lambda \sqrt{\mu} + \mu + c^2 \l \xi^2 + \alpha^2 \r}.
\label{fourierlaplace un mezzo bid}
\end{align}
In the previous calculation we use the Laplace transform of $\mathpzc{l}_\frac{1}{2} (x, t)$ obtained in \eqref{tlaplace l un mezzo}. The proof is complete since \eqref{fourierlaplace un mezzo bid},  coincides with \eqref{fourier laplace proc plan comp} and with the Fourier-Laplace transform of \eqref{timw frac bid}.
\end{proof}
\end{te}
\begin{os} \normalfont
Since for the first passage time $\tau_{\frac{s}{\sqrt{2}}} = \inf \ll z:B(z) = \frac{s}{\sqrt{2}} \rr$ of a Brownian motion through level $\frac{s}{\sqrt{2}}$ we have that
\begin{equation} 
	\int_0^\infty e^{-\mu t} \Pr \ll \tau_{\frac{s}{\sqrt{2}}} \in dt \rr \,  = \, e^{-s\sqrt{\mu}},
	\label{}
\end{equation}
and
\begin{equation} 
	\int_0^\infty e^{-\mu t} \frac{\partial^{\frac{1}{2}}}{\partial t^{\frac{1}{2}}} p_{|B|} (s, t) \, dt \, = \, e^{-s\sqrt{\mu}}
	\label{}
\end{equation}
we can write
\begin{align}
	&  \int_0^\infty \mathfrak{r} (x, y, s) \frac{\partial^{\frac{1}{2}}}{\partial t^{\frac{1}{2}}} p_{|B|} (s, t) \, ds \, = \, \int_0^\infty \mathfrak{r} (x, y, s) \frac{s}{\sqrt{2}} \frac{e^{-\frac{s^2}{4t}}}{\sqrt{2\pi t^3}} \, ds \notag \\
	 = \, & \int_0^\infty \frac{\partial}{\partial s} \mathfrak{r} (x, y, s) \frac{e^{-\frac{s^2}{4t}}}{\sqrt{\pi t}} ds \, = \, \int_0^\infty \frac{\partial}{\partial s} \mathfrak{r}(x, y, s) \, p_{|B|} (s, t) \, ds.
	\label{due asterischi}
\end{align}
This representation of the  second term of \eqref{processo finito composto} is extremely interesting because by integrating \eqref{due asterischi} in $C_{ct}$ we get
\begin{equation} 
	\int_0^\infty \frac{\partial}{\partial s} (1-e^{-2\lambda s}) p_{|B|} (s, t) \, ds \, = \,  2\lambda \int_0^\infty e^{-2\lambda s} p_{|B|} (s, t) \, ds
	\label{487}
\end{equation}
and yields the missing probability of the first term of \eqref{processo finito composto}.
\end{os}

\begin{os} \normalfont
We check that the law
\begin{equation} 
	\mathfrak{q} (x, y, t) \, = \, \int_0^\infty \mathfrak{r} (x, y, s) \left[ p_{|B|} (s, t) + \frac{1}{2\lambda} \frac{\partial^{\frac{1}{2}}}{\partial t^{\frac{1}{2}}} p_{|B|} (s, t) \right] \, ds
	\label{}
\end{equation}
integrates to unity. By taking the $t$-Laplace transform, the integral with respect to $(x, y)$ becomes
\begin{align}
	&  \iint_{C_{ct}} dx \, dy \int_0^\infty dt \, e^{-\mu t}  \, \mathfrak{q} (x, y, t) \notag \\
	 = \, & \int_0^\infty \l 1-e^{-2\lambda s} \r \left[ \int_0^\infty e^{-\mu t} \l p_{|B|} (s, t) + \frac{1}{2\lambda} \frac{\partial^{\frac{1}{2}}}{\partial t^{\frac{1}{2}}} p_{|B|} (s, t) \r dt \right] ds \notag \\
	 = \, & \int_0^\infty \l 1-e^{-2\lambda s} \r \left[ \frac{e^{-s\sqrt{\mu }}}{\sqrt{\mu }} + \frac{e^{-s\sqrt{\mu }}}{2\lambda} \right] \, ds \notag \\
	 = \, & \l \frac{1}{\sqrt{\mu}} + \frac{1}{2\lambda} \r \left[ \int_0^\infty e^{-s\sqrt{\mu }} ds - \int_0^\infty e^{-s \l 2\lambda + \sqrt{\mu} \r} ds \right] \notag \\
	 = \, & \frac{2\lambda + \sqrt{\mu }}{2\lambda \sqrt{\mu}} \l \frac{1}{\sqrt{\mu}} - \frac{1}{2\lambda + \sqrt{\mu}} \r \, = \, \frac{1}{\mu} \, = \, \int_0^\infty e^{-\mu t} dt.
	\label{}
\end{align}
The same check can be done directly by taking into account formulas \eqref{due asterischi} and \eqref{487}.
\end{os}
Relationships similar to $B \l c^2 \mathpzc{L}^{\frac{1}{2}} (t) \r \, \stackrel{\textrm{law}}{=} \, T \l \lll B(t) \rrr \r$, $t>0$, and the analogous one in the plane, cannot be established in spaces of dimension $n \geq 3$, because random motions governed by telegraph equations in such spaces have not been constructed. Random flights in $\mathbb{R}^n$ have been studied (\citet{orsdegr}) but their distributions are not related to higher-dimensional telegraph equations.

\end{document}